\documentclass[12pt]{amsart}
\usepackage{amsmath, amssymb, latexsym, amsthm, mathrsfs, color,  geometry}
\usepackage[all]{xy}
\newcommand{\ed}{

\subsection*{Acknowledgments}
Parts of the work reported here were carried out during two visits of the second named author
at the Kurt G\"odel Research Center.
These visits were partially supported by FWF Grant M1244-N13 and
by the Research Networking Programme \emph{New Frontiers of Infinity} (INFTY),
funded by the European Science Foundation.

The second named author thanks the third named author for his kind hospitality,
and the Kurt G\"odel Research Center Director, researchers and staff for the excellent academic and friendly atmosphere.

\end{document}}

      \newenvironment{changemargin}[2]{\begin{list}{}{
         \setlength{\topsep}{0pt}\setlength{\leftmargin}{0pt}
         \setlength{\rightmargin}{0pt}
         \setlength{\listparindent}{\parindent}
         \setlength{\itemindent}{\parindent}
         \setlength{\parsep}{0pt plus 1pt}
         \addtolength{\leftmargin}{#1}\addtolength{\rightmargin}{#2}
         }\item }{\end{list}}

\newcommand{\nc}{\newcommand}
\newcommand{\two}{\{0,1\}}
\nc{\productive}[2]{\bigl(#1,\allowbreak #2\bigr)^\x}
\nc{\set}[2]{{\left\{#1 \;:\; #2\right\}}}
\nc{\seq}[2]{{\la #1 \;:\; #2\ra}}
\nc{\cube}{(\Cantor)^\N}
\nc{\Match}{\op{Match}}
\nc{\concat}[1]{\hat{\phantom{a}}\langle #1\rangle}
\nc{\poset}{\mathbb{P}}
\nc{\fn}[1]{{\op{Fn}(#1\times\w,2)}}
\nc{\linadd}{\op{linadd}}
\nc{\nonprod}{\non^\x}
\nc{\Ga}{\Gamma}
\nc{\Om}{\Omega}
\nc{\alephes}{{\aleph_0}}
\nc{\my}[1]{{\color{red} #1}}
\nc{\Cp}{\op{C}_p}
\nc{\Bp}{\op{B}_p}
\nc{\Pa}[8]{\bibitem{#1} {#2}, \emph{#3}, {#4} \textbf{#5} ({#6}), {#7}--{#8}.}
\nc{\tPa}[5]{\bibitem{#1} {#2}, \emph{#3}, {#4}, to appear.}
\nc{\sPa}[4]{\bibitem{#1} {#2}, \emph{#3}, {#4}, submitted.}
\nc{\Bc}[9]{\bibitem{#1} {#2}, \emph{#3}, in: \textbf{#4} (#5), #6 #7, #8--#9.}
\nc{\fD}{\mathfrak{D}}
\nc{\fX}{\mathfrak{X}}
\nc{\Onbd}{\Op_{\mathrm{nbd}}} %{\Op_{\mathsf{nbd}}}
\nc{\Omnb}{\Om_{\mathrm{nbd}}} %{\Om_{\mathsf{nbd}}}
\nc{\od}{\mathfrak{od}}
\nc{\Setting}[7]{\xymatrix@R=4pt@C=7pt{#1\ar@{-}[r]&#2\ar@{-}[r]&#3\\&#4\ar@{-}[u]\\
#5\ar@{-}[uu]\ar@{-}[r] & #6\ar@{-}[u]\ar@{-}[r] & #7\ar@{-}[uu]}}
\nc{\mx}[1]{\begin{matrix}#1\end{matrix}}
\nc{\plim}{p\txt{-}\lim}
\nc{\Bgp}{{\Z^\N}}
\nc{\Cgp}{{{\Z_2}^\N}}
\nc{\Cite}[1]{\textbf{[#1]}}
\nc{\Next}[1]{{#1^+}}
\nc{\Fr}{\mathit{F\!r}}
\nc{\intvl}[2]{{[#1(#2),\allowbreak #1(#2\!+\!1))}}
\nc{\Bdd}{\mathbf{B}}
\nc{\Ax}{\mathsf{Ax}}
\nc{\Dfin}{\mathfrak{D}_\mathrm{fin}}
\nc{\grbl}{{\mbox{\textit{\tiny gp}}}}
\nc{\bbP}{\mathbb{P}}
\nc{\BOfat}{\B_{\Om_{\mathrm{fat}}}}%\B_{\mathrm{fat}}}
\nc{\Bgood}{\B_{\mathrm{good}}}
\nc{\compactN}{\cl{\mathbb{N}}}
\nc{\blocks}[2]{\op{cl}_{#2}(#1)}
\nc{\blocksplus}[2]{\op{cl}^+_{#2}(#1)}
\nc{\arx}[1]{\texttt{http://arxiv.org/math/#1}}
\nc{\bq}{\begin{quote}}
\nc{\eq}{\end{quote}}
\nc{\cl}[1]{\overline{#1}}
\nc{\CH}{the Continuum Hypothesis}
\nc{\MA}{Martin's Axiom}
\nc{\Bfat}{\B_\mathrm{fat}}
\nc{\inv}{^{-1}}
\nc{\Cantor}{{\two^\N}}%{{2^\w}}
\nc{\bP}{\mathbf{P}}
\nc{\bof}{\op{\fb}}
\nc{\bofF}{\bof(\cF)}
\nc{\sr}[3]{{\txt{#1\\#3}}}%{{\txt{#1\\#2\\#3}}}
\nc{\gp}{\binom{\Om}{\Ga}}
\nc{\gpsmall}{\mbox{$\gp$}}
\nc{\gig}{\gimel}%{\gimel\Ga}
\nc{\gns}{\sone(\Om,\gig)}
\nc{\nsr}[2]{#1}
\nc{\N}{\mathbb{N}}%{\w}
\nc{\NN}{{\N^{\N}}}
\nc{\ZN}{{\Z^{\N}}}
\nc{\NNup}{{\N^{\uparrow\N}}}
\nc{\PN}{{P(\N)}}
\nc{\roth}{{[\N]^{\mbox{\tiny $\infty$}}}} %{{[\w]^{\w}}}
\nc{\Fin}{[\N]^{\mbox{\tiny $<\!\infty$}}} %{{[\w]^{<\w}}}%{{[\N]^{<\aleph_0}}}
\nc{\ici}{{[\w]^{(\w,\w)}}}%{{[\N]^{(\aleph_0,\aleph_0)}}}
\nc{\Inc}{{\compactN^{\uparrow\N}}}
\nc{\powInc}[1]{{\big(\Inc\big)^{#1}}}
\nc{\powFin}[1]{{\big(\Fin\big)^{#1}}}
\nc{\powPN}[1]{{\big(\PN\big)^{#1}}}
\nc{\NcompactN}{{\compactN^\N}}
\nc{\setseq}[1]{\{#1 : n\in\N\}}
\nc{\sseq}[1]{\{#1 : n\in\N\}}
\nc{\Uarrow}{\smash{\big\uparrow}}
\nc{\LE}{\preccurlyeq}
\nc{\GE}{\succcurlyeq}
\nc{\op}{\operatorname}
\nc{\im}{\op{im}}
\nc{\Span}{\op{span}}
\nc{\maxfin}{\op{maxfin}}
\nc{\ran}{\op{range}}
\nc{\iso}{\cong}
\nc{\Madd}{{\M}^\star}
\nc{\cI}{\mathcal{I}}
\nc{\cJ}{\mathcal{J}}
\nc{\scrA}{\mathscr{A}}
\nc{\scrB}{\mathscr{B}}
\nc{\scrC}{\mathscr{C}}
\nc{\scrD}{\mathscr{D}}
\nc{\scrF}{\mathscr{F}}
\nc{\scrK}{\mathscr{K}}
\nc{\A}{\forall}
\nc{\B}{\mathrm{B}}
\nc{\cB}{\mathcal{B}}
\nc{\bB}{\mathbf{B}}
\nc{\BG}{\B_\Ga}
\nc{\BL}{\B_\Lambda}
\nc{\BT}{\B_\Tau}
\nc{\BTstar}{\B_{\Tau^*}}
\nc{\BO}{\B_\Om}
\nc{\DO}{\cD_\Om}
\nc{\KO}{\cK_\Om}
\nc{\CG}{C_\Ga}
\nc{\CL}{C_\Lambda}
\nc{\CT}{C_\Tau}
\nc{\CTstar}{C_{\Tau^*}}
\nc{\CO}{C_\Om}
\nc{\COgp}{C_{\Om^{\grbl}}}
\nc{\CLgp}{C_{\Lambda^{\grbl}}}
\nc{\BOgp}{\B_{\Om}^{\grbl}}
\nc{\BLgp}{\B_{\Lambda^{\grbl}}}
\nc{\sfC}{\mathsf{C}}
\nc{\sfD}{\mathsf{D}}
\nc{\bD}{\mathbf{D}}
\nc{\Tau}{\mathrm{T}}
\nc{\cA}{\mathcal{A}}
\nc{\cK}{\mathcal{K}}
\nc{\cD}{\mathcal{D}}
\nc{\cF}{\mathcal{F}}
\nc{\cS}{\mathcal{S}}
\nc{\cG}{\mathcal{G}}
\nc{\cY}{\mathcal{Y}}
\nc{\J}{\mathcal{J}}
\nc{\cL}{\mathcal{L}}
\nc{\cM}{\mathcal{M}}
\nc{\cN}{\mathcal{N}}
\nc{\cO}{\mathcal{O}}
\nc{\Op}{\mathrm{O}}
\nc{\cP}{\mathcal{P}}
\nc{\Q}{\mathbb{Q}}
\nc{\R}{\mathbb{R}}
\nc{\cU}{\mathcal{U}}
\nc{\Un}{\bigcup}
\nc{\cV}{\mathcal{V}}
\nc{\cW}{\mathcal{W}}
\nc{\Z}{{\mathbb Z}}
\nc{\Impl}{\Rightarrow}
\long\def\forget#1\forgotten{}
\nc{\ft}{\mathfrak{t}}
\nc{\fb}{\mathfrak{b}}
\nc{\fc}{\mathfrak{c}}
\nc{\fd}{\mathfrak{d}}
\nc{\fg}{\mathfrak{g}}
\nc{\oo}{\infty}
\nc{\fr}{\mathfrak{r}}
\nc{\fu}{\mathfrak{u}}
\nc{\fh}{\mathfrak{h}}
\nc{\fp}{\mathfrak{p}}
\nc{\fj}{\mathfrak{j}}
\nc{\fs}{\mathfrak{s}}
\nc{\w}{\omega}
\nc{\x}{\times}
\nc{\Iff}{\Leftrightarrow}
\newcommand\comp{^{\text{\tt c}}}
\nc{\nin}{\notin}
\nc{\cat}{\hat{\ }}
\nc{\sub}{\subseteq}
\nc{\spst}{\supseteq}
\nc{\sm}{\setminus}
\nc{\as}{\subseteq^*}%{\let\proclaim\relax}
\nc{\rest}{\restriction}

\nc{\la}{\langle}
\nc{\ra}{\rangle}
\nc{\E}{\exists}
\nc{\dom}{\op{dom}}
\nc{\cov}{\op{cov}}
\nc{\add}{\op{add}}
\nc{\cof}{\op{cof}}
\nc{\cf}{\op{cf}}
\nc{\non}{\op{non}}
\nc{\unif}{\op{non}}
\nc{\COV}{\op{COV}}
\nc{\ADD}{\op{ADD}}
\nc{\COF}{\op{COF}}
\nc{\NON}{\op{NON}}
\nc{\impl}{\to}
\nc{\Lp}{\mathcal{L_\p}}
\nc{\Wlog}{without loss of generality}
\newtheorem{thm}{Theorem}[section]
\nc{\bthm}{\begin{thm}} \nc{\ethm}{\end{thm}}
\newtheorem{prop}[thm]{Proposition}
\nc{\bprp}{\begin{prop}} \nc{\eprp}{\end{prop}}
\newtheorem{fact}[thm]{Fact}
\nc{\bfct}{\begin{fact}} \nc{\efct}{\end{fact}}
\newtheorem{prob}[thm]{Problem}
\nc{\bprb}{\begin{prob}} \nc{\eprb}{\end{prob}}
\newtheorem{lem}[thm]{Lemma}
\nc{\blem}{\begin{lem}} \nc{\elem}{\end{lem}}
\newtheorem{claim}[thm]{Claim}
\nc{\bclm}{\begin{claim}} \nc{\eclm}{\end{claim}}
\newtheorem{cor}[thm]{Corollary}
\nc{\bcor}{\begin{cor}} \nc{\ecor}{\end{cor}}
\newtheorem{conj}[thm]{Conjecture}
\nc{\bcnj}{\begin{conj}} \nc{\ecnj}{\end{conj}}
\theoremstyle{definition}
\newtheorem{defn}[thm]{Definition}
\nc{\bdfn}{\begin{defn}} \nc{\edfn}{\end{defn}}
\theoremstyle{remark}
\newtheorem{rem}[thm]{Remark}
\nc{\brem}{\begin{rem}} \nc{\erem}{\end{rem}}
\newtheorem{cnv}[thm]{Convention}
\nc{\bcnv}{\begin{cnv}} \nc{\ecnv}{\end{cnv}}
\newtheorem{exam}[thm]{Example}
\nc{\bexm}{\begin{exam}} \nc{\eexm}{\end{exam}}
\nc{\bpf}{\begin{proof}} \nc{\epf}{\end{proof}}
\nc{\be}{\begin{enumerate}}
\nc{\ee}{\end{enumerate}}
\nc{\bi}{\begin{itemize}}
\nc{\itm}{\item}
\nc{\ei}{\end{itemize}}
\nc{\Subsection}[1]{\goodbreak\subsection*{#1}}%\ \par}

\forget
\forgotten

\nc{\sone}{\mathsf{S}_1}
\nc{\sfin}{\mathsf{S}_\mathrm{fin}}
\nc{\ufin}{\mathsf{U}_\mathrm{fin}}
\nc{\Split}{\mathsf{Split}}

\nc{\gone}{\mathsf{G}_1}    \nc{\gfin}{\mathsf{G}_\mathrm{fin}}

\title[Coverings and products]{Selective covering properties of product spaces}
\author[Arnold Miller]{Arnold W. Miller}
\address[Miller]{Department of Mathematics, University of Wisconsin-Madison, Van Vleck Hall 480 Lincoln Drive, Madison, Wisconsin 53706-1388, USA}
\email{miller@math.wisc.edu}
\urladdr{http://www.math.wisc.edu/\~{}miller/}
\author{Boaz Tsaban}
\address[Tsaban]{Department of Mathematics, Bar-Ilan University, Ramat Gan 52900, Israel}
\email{tsaban@math.biu.ac.il}
\urladdr{http://www.cs.biu.ac.il/\~{}tsaban}
\author{Lyubomyr Zdomskyy}
\address[Zdomskyy]{Kurt G\"odel Research Center for Mathematical Logic, University of Vienna,
W\"ahringer Str.\ 25, 1090 Vienna, Austria}
\email{lzdomsky@logic.univie.ac.at}
\urladdr{http://www.logic.univie.ac.at/\~{}lzdomsky/}

\keywords{%
Gerlits--Nagy $\gamma$ property, 
Menger property, Hurewicz property, Rothberger property, Gerlits-Nagy $(*)$ property,
productively Lindel\"of,
product theory, selection principles, special sets of real numbers. 
}

\subjclass{%
26A03, %Foundations: limits and generalizations, elementary topology of the line
03E17, %Cardinal characteristics of the continuum
03E75 %Applications of set theory
}

\begin{document}

\begin{abstract}
We study the preservation of selective covering properties, including classic ones introduced by
Menger, Hurewicz, Rothberger, Gerlits and Nagy,
and others, under products with some major families of concentrated sets of reals.

Our methods include the projection method introduced by the authors in an earlier work, as well as several new
methods. Some special consequences of our main results are (definitions provided in the paper):
\be
\item Every product of a concentrated space with a Hurewicz $\sone(\Ga,\Op)$ space satisfies $\sone(\Ga,\Op)$. On the other hand,
assuming \CH{}, for each Sierpi\'nski set $S$ there is a Luzin set $L$ such that $L\x S$ can be mapped onto the real line
by a Borel function.
\item Assuming Semifilter Trichotomy, every concentrated space is productively Menger and productively Rothberger.
\item Every scale set is productively Hurewicz, productively Menger, productively Scheepers, and productively Gerlits--Nagy.
\item Assuming $\fd=\aleph_1$, every productively Lindel\"of space is productively Hurewicz, productively Menger, and productively Scheepers.
\ee
A notorious open problem asks whether the additivity of Rothberger's property may be strictly greater than $\add(\cN)$, the additivity of the ideal of Lebesgue-null sets of reals.
We obtain a positive answer, modulo the consistency of Semifilter Trichotomy with $\add(\cN)<\cov(\cM)$.

Our results improve upon and unify a number of results, established earlier by many authors.
 \end{abstract}

\maketitle

\section{Introduction}

All topological spaces in this paper are assumed, without further mention, to be Tychonoff.
Since the results presented here are new even in the case where the spaces are subsets of the real line,
readers who wish to do so may assume throughout that we deal with sets of real numbers.

We study selective covering properties of products of topological spaces.
Our results, that answer questions concerning classic covering properties,
are best perceived in the modern framework of selection principles,
to which we provide here a brief introduction.\footnote{This introduction
is adopted from \cite{LinSAdd}. Extended introductions
to this field are available in \cite{KocSurv, LecceSurvey, ict}.}
This framework was introduced by Scheepers in \cite{coc1}
to study, in a uniform manner, a variety of properties introduced
in different mathematical disciplines, since the early
1920's, by Menger, Hurewicz, Rothberger, Gerlits and Nagy, and others.

Let $X$ be a topological space. We say that $\cU$ is a \emph{cover}
of $X$ if $X=\Un\cU$, but $X\nin\cU$.
Often, $X$ is considered as a subspace of another space $Y$,
and in this case we always consider covers of $X$ by subsets of $Y$,
and require instead that no member of the cover contains $X$.
Let $\Op(X)$ be the family of open covers of $X$.
Define the following subfamilies of $\Op(X)$:
$\cU\in\Om(X)$ if each finite subset of $X$ is contained in some member of $\cU$.
$\cU\in\Ga(X)$ if $\cU$ is infinite, and each element of $X$ is contained in all but
finitely many members of $\cU$.

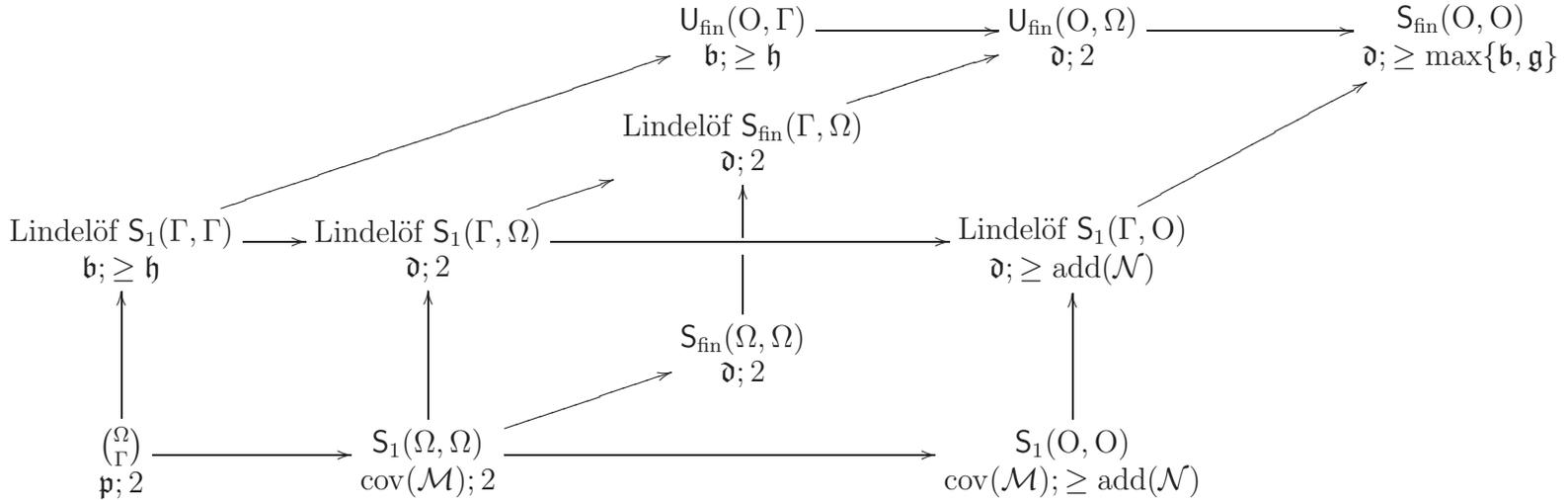
\begin{figure}[!htp]
\begin{changemargin}{-4cm}{-3cm}
\begin{center}
{%\scriptsize
$\xymatrix@R=8pt{
%1
&
&
& \sr{$\ufin(\Op,\Ga)$}{Hurewicz}{$\fb; \ge\fh$}\ar[r]
& \sr{$\ufin(\Op,\Om)$}{}{$\fd; 2$}\ar[rr]
& & \sr{$\sfin(\Op,\Op)$}{Menger}{$\fd; \ge\max\{\fb,\fg\}$}
\\
%2
&
&
& \sr{Lindel\"of $\sfin(\Ga,\Om)$}{}{$\fd; 2$}\ar[ur]
\\
%3
& \sr{Lindel\"of $\sone(\Ga,\Ga)$}{}{$\fb; \ge \fh$}\ar[r]\ar[uurr]
& \sr{Lindel\"of $\sone(\Ga,\Om)$}{}{$\fd; 2$}\ar[rr]\ar[ur]
& & \sr{Lindel\"of $\sone(\Ga,\Op)$}{}{$\fd; \ge\add(\cN)$}\ar[uurr]
\\
%4
&
&
& \sr{$\sfin(\Om,\Om)$}{}{$\fd; 2$}\ar'[u][uu]
\\
%5
& \sr{$\gp$}{Gerlits--Nagy}{$\fp; 2$}\ar[r]\ar[uu]
& \sr{$\sone(\Om,\Om)$}{}{$\cov(\cM); 2$}\ar[uu]\ar[rr]\ar[ur]
& & \sr{$\sone(\Op,\Op)$}{Rothberger}{$\cov(\cM); \ge\add(\cN)$}\ar[uu]
}$
}
\caption{The Scheepers Diagram}\label{SchDiag}
\end{center}
\end{changemargin}
\end{figure}

Some of the following statements may hold for families $\scrA$ and $\scrB$ of covers of $X$.
\begin{description}
\item[$\binom{\scrA}{\scrB}$] Each member of $\scrA$ contains a member of $\scrB$.
\item[$\sone(\scrA,\scrB)$] For each sequence $\seq{\cU_n\in\scrA}{n\in\N}$, there is a selection
$\seq{U_n\in\cU_n}{n\in\N}$ such that $\sseq{U_n}\in\scrB$.
\item[$\sfin(\scrA,\scrB)$] For each sequence $\seq{\cU_n\in\scrA}{n\in\N}$, there is a selection
of finite sets $\seq{\cF_n\sub\cU_n}{n\in\N}$ such that $\Un_n\cF_n\in\scrB$.
\item[$\ufin(\scrA,\scrB)$] For each sequence $\seq{\cU_n\in\scrA}{n\in\N}$, where no $\cU_n$ contains a finite
subcover, there is a selection
of finite sets $\seq{\cF_n\sub\cU_n}{n\in\N}$ such that $\sseq{\Un\cF_n}\in\scrB$.
\end{description}
We say, e.g., that $X$ satisfies $\sone(\Op,\Op)$ if the statement $\sone(\Op(X),\Op(X))$ holds.
This way, $\sone(\Op,\Op)$ is a property (or a class) of topological spaces, and similarly for all other statements
and families of covers.
In the realm of Lindel\"of spaces,\footnote{Indeed, all properties in the Scheepers Diagram \ref{SchDiag},
except for those having $\Gamma$ in the first argument, imply being Lindel\"of.}
each nontrivial property among these properties, where $\scrA,\scrB$ range over $\Op,\Om,\Ga$,
is equivalent to one in Figure \ref{SchDiag}  \cite{coc1, coc2}.
In this diagram, an arrow denotes implication.

The extremal properties in this diagram are classic and were introduced by Menger ($\sfin(\Op,\Op)$, under a difference guise),
Hurewicz ($\ufin(\Op,\Ga)$), Rothberger ($\sone(\Op,\Op)$), and Gerlits and Nagy ($\binom{\Om}{\Ga}$). The other ones
were introduced and studied more recently, by many authors.

In this diagram, we indicate below each class $P$ its \emph{critical cardinality} $\non(P)$,
the minimal cardinality of a space not in the class, and its (provable) \emph{additivity number} $\add(P)$,
the minimal number of spaces possessing this property, whose union does \emph{not} have this property.
These cardinals are all combinatorial cardinal characteristics of the continuum,
details about which are available in \cite{BlassHBK}. Here,
$\cM,\cN$ are the families of meager (i.e., Baire first category) sets in $\R$ and Lebesgue null sets in $\R$,
respectively.
In cases where only lower bounds on $\add(P)$ are given, $\cf(\non(P))$ is an upper bound.

On occasions, we will also consider the classes of covers
$\B$, $\BO$ and $\BG$, defined as $\Op$, $\Om$ and $\Ga$ were
defined, replacing \emph{open cover} by \emph{countable Borel cover}.
The properties thus obtained have rich history of their own \cite{CBC}, and for Lindel\"of spaces,
the Borel variants of the studied properties are (usually, strictly) stronger than the open ones \cite{CBC}.

Many additional---classic and new---properties were or can be studied in relation to the Scheepers Diagram.
Some examples of this kind are provided in the present paper.

\medskip

The following definition and observation are useful.

\bdfn
Let $P$ be a property (or class) or topological spaces.
A topological space $X$ is \emph{productively $P$} if $X\x Y$ has the property $P$ for each $Y$ satisfying $P$.
$P^\uparrow$ is the property of having all finite powers satisfying $P$.
\edfn

In this notation, $\sone(\Op,\Op)^\uparrow=\sone(\Om,\Om)$ \cite{Sakai88} and
$\sfin(\Op,\Op)^\uparrow=\sfin(\Om,\Om)$ \cite{coc2}.
If $X$ is productively $P$ and the singleton space satisfies $P$, then
$X$ satisfies $P$. Moreover, we have the following.

\blem\label{prodispow}
Let $X$ be a productively $P$ topological space. Then:
\be
\itm Every finite power of $X$ is productively $P$.
\itm $X$ is productively $P^\uparrow$.
\itm Every finite power of $X$ is productively $P^\uparrow$.
\ee
\elem
\bpf
(1) By induction on the power of $X$, $X^k\x Y$ has the property $P$ if $Y$ has it.

(2) Let $Y$ be in $P^\uparrow$. For each $k$, $(X\x Y)^k\iso
 X^k\x Y^k$. Apply (1).

(3) By (1) and (2).
\epf

In particular, if $X$ is productively $\sone(\Op,\Op)$, then it is also productively $\sone(\Om,\Om)$,
and similarly for $\sfin$. Several additional properties in the literature are characterized by
having the form $P^\uparrow$ for a property $P$ in the Scheepers Diagram, and the same comment
applies.

\section{Concentrated spaces and $\sone(\Ga,\Ga)$}

Let $\kappa$ be an uncountable cardinal.
Following Besicovich \cite{Besicovitch1, Besicovitch2}, we say that a topological space $X$ is
$\kappa$-\emph{concentrated} if there is a countable set $D\sub X$ such that
$|X\sm U|<\kappa$ for every open set $U\spst D$.
Several major examples of families of concentrated spaces will be considered later.

A special case of Theorem 11(3) in Babinkostova--Scheepers \cite{BabSch} is that for each concentrated metric space $C$,
if $Y$ satisfies $\ufin(\Op,\Ga)$ and $\sone(\Op,\Op)$, then $C\x Y$ satisfies $\sone(\Op,\Op)$.
Theorem 3.1 in the more recent paper \cite{AddGN} implies, in particular, that it suffices to assume that $C$ is a
$\cov(\cM)$-concentrated space.
Our first observation is that the methods of the paper \cite{AddGN} imply a similar result for $\sone(\Ga,\Ga)$.
The proof given here is slightly more general than the one that may be extracted from \cite{AddGN}.

\blem\label{addgggo}
Let a topological space $X$ be a union of less than $\cov(\cM)$ many $\sone(\Ga,\Ga)$ spaces.
Then $X$ satisfies $\sone(\Ga,\Op)$.
\elem
\bpf
Write $X=\Un_{\alpha<\kappa}X_\alpha$, with $\kappa<\cov(\cM)$ and each $X_\alpha$ satisfying $\sone(\Ga,\Ga)$.
Let $\cU_n\in\Ga(X)$ for all $n$. As $\Ga(X)$ is closed taking infinite subsets, we may assume that we can enumerate
$\cU_n=\set{U^n_m}{m\in\N}$ for all $n$.

For each $\alpha<\kappa$, let $f_\alpha\in\NN$ be such that $\set{U^n_{f_\alpha(n)}}{n\in\N}\in\Ga(X_\alpha)$.
As $\kappa<\cov(\cM)$, there is $f\in\NN$ such that, for each $\alpha<\kappa$, $f(n)=f_\alpha(n)$ for infinitely many $n$.
Then $\set{U^n_{f(n)}}{n\in\N}$ is a cover of $X$.
\epf

The method used in the following proof was introduced in \cite{AddGN}. Since this method is used frequently in the present
paper, we name it \emph{the projection method}.

\bthm\label{concgggo}
Let $C$ be a $\cov(\cM)$-concentrated space. For each Lindel\"of $\sone(\Ga,\Ga)$ space $Y$,
$X\x Y$ satisfies $\sone(\Ga,\Op)$.
\ethm
\bpf
Let $C$ be $\cov(\cM)$-concentrated on some countable set $D\sub C$.
Let $Y$ be a Lindel\"of $\sone(\Ga,\Ga)$ space.
Let $K$ be a compact space containing $C$ as a subspace.
For each $n$, let $\cU_n\in\Ga(C\x Y)$, where the elements of $\cU_n$ are open in $K\x Y$.

As $D$ is countable, $D\x Y$ satisfies $\sone(\Ga,\Op)$ (Lemma \ref{addgggo}).
Pick $U_n\in\cU_n$, $n\in\N$, such that $D\x Y\sub U:=\Un_nU_n$.

The Hurewicz property $\ufin(\Op,\Ga)$ is preserved by products with compact spaces, moving to
closed subspaces,
and continuous images \cite{coc2}.
Since $Y$ satisfies $\ufin(\Op,\Ga)$ and $K$ is compact, $K\x Y$ satisfies $\ufin(\Op,\Ga)$.
Thus, so does $K\x Y\sm U$. It follows that the projection $H$ of $(K\x Y)\sm U$
on the first coordinate, satisfies $\ufin(\Op,\Ga)$.
Note that
$$(K\sm H)\x Y\sub U.$$
The argument in the proof of \cite[Theorem 5.7]{coc2} generalizes to
arbitrary (Tychonoff) spaces, to show that for $H,F$ disjoint subspaces of a space $K$
with $H$ $\ufin(\Op,\Ga)$, and $F$ F$_\sigma$, there is a G$_\delta$ set $G\sub K$
such that $G\spst F$ and $H\cap G=\emptyset$.

Let $G$ be a G$_\delta$ subset of $K$ such that
$D\sub G$ and $H\cap G=\emptyset$. As $C$ is $\cov(\cM)$-concentrated on
$D$, $C\sm G$ is a countable increasing union of sets of cardinality $<\cov(\cM)$.
By Lemma \ref{addgggo} and the fact that $\sone(\Ga,\Op)$ is preserved under countable
unions (or, alternatively, as $\cov(\cM)$ has uncountable cofinality),
$(C\sm G)\x Y$ satisfies $\sone(\Ga,\Op)$. Take $V_n\in\cU_n$,
$n\in\N$, such that $(C\sm G)\x Y\sub \Un_n V_n$.
Then
$$C\x Y\sub\Un_{n\in\N}(U_n\cup V_n).$$
We have picked two sets (instead of one) from each cover $\cU_n$, but this is fine (e.g., \cite[Appendix A]{MHP}).
\epf

The methods of \cite{AddGN} also imply the following, more general result. Since the proof is similar
to that of Theorem 3.3 in \cite{AddGN} and we are not going to use this result here, we omit the proof.

\bdfn
Let $\kappa$ be an infinite cardinal number.
Let $\scrC_0(\kappa)$ be the family of spaces of cardinality $<\kappa$. For successor ordinals $\alpha+1$,
let $C\in\scrC_{\alpha+1}(\kappa)$ if:
\be
\itm either there is a countable $D\sub C$ with $C\sm U\in\scrC_\alpha(\kappa)$ for all open $U\spst D$;
\itm or $C$ is a union of less than $\cf(\kappa)$ members of $\scrC_\alpha(\kappa)$.
\ee
For limit ordinals $\alpha$, let $\scrC_\alpha(\kappa)=\Un_{\beta<\alpha}\scrC_\beta(\kappa)$.
\edfn

By the Scheepers Diagram,  $\add(\cN)\le\add(\sone(\Ga,\Op))$.

\bthm\label{conc3}
The product of each member of $\scrC_{\add(\cN)}(\cov(\cM))$
with every Lindel\"of $\sone(\Ga,\Ga)$ space satisfies $\sone(\Ga,\Op)$.\qed
\ethm

%The proof of the following theorem, and that of the forthcoming Theorem \ref{conc4},
%use that $\add(\cN)\le\add(\sone(\Ga,\Op))$ \cite{AddGN}.

\bdfn
Let $P,Q$ be classes of spaces, each containing all one-element spaces and closed under
homeomorphic images.
$\productive{P}{Q}$ is the class of all spaces $X$ such that, for each $Y$ in $P$,
$X\x Y$ is in $Q$.
%$\productive{P}{P}$ is denoted $P^\x$.
\edfn

By Lemma \ref{addgggo}, $\cov(\cM)\le\non(\productive{\sone(\Ga,\Ga)}{\sone(\Ga,\Op)})$.
Theorem \ref{conc3} holds, more generally, for $\scrC_{\add(\sone(\Ga,\Op))}(\non(\productive{\mbox{Lindel\"of }\sone(\Ga,\Ga)}{\sone(\Ga,\Op)}))$.

\bprb
Is $\non(\productive{\mbox{Lindel\"of }\sone(\Ga,\Ga)}{\sone(\Ga,\Op)})=\fd$?
\eprb

\section{Concentrated sets and the conjunction of $\ufin(\Op,\Ga)$ and $\sone(\Ga,\Op)$}

In this section, we consider the conjunction of $\ufin(\Op,\Ga)$ and $\sone(\Ga,\Op)$.
This class is larger than Lindel\"of $\sone(\Ga,\Ga)$.
The definition of \emph{$\fb$-scale set} is given in Section \ref{sec:scales}. For the present purpose, it suffices
to know their following properties (cf.\ \cite{MHP}): $\fb$-scale sets are subspaces of $\R$, of cardinality $\fb$,
that can be constructed outright in ZFC.
They are $\fb$-concentrated, and as such satisfy  $\sone(\Ga,\Op)$, and they satisfy  $\ufin(\Op,\Ga)$.
The following results are known.

\bthm
\mbox{}
\be
\itm Every $\fb$-scale set satisfies  $\ufin(\Op,\Ga)$ and $\sone(\Ga,\Op)$ \cite{ideals} (cf.\ \cite{MHP}).
\itm Consistently, no set of reals of cardinality $\fb$ satisfies $\sone(\Ga,\Ga)$ \cite{BBC}.
\itm The Continuum Hypothesis implies that there is a $\fb$-scale set not satisfying $\sone(\Ga,\Ga)$ \cite{HurPerf}.
\ee
\ethm

We will show that the conjunction of $\ufin(\Op,\Ga)$ and $\sone(\Ga,\Op)$
can be expressed as a standard selective property.
A countable cover $\cU$ of a space $X$ is in $\gimel(\Gamma)$ \cite{GlCovs}
if for each (equivalently, some)
bijective enumeration $\cU=\sseq{U_n}$, there is an increasing $h\in\NN$ such that, for each $x\in X$,
$$x\in \Un_{k=h(n)}^{h(n+1)-1}U_k$$
for all but finitely many $n$.
In \cite{coc7} it is shown that $\ufin(\Op,\Ga)\cap \sone(\Op,\Op)=\sone(\Om,\gimel(\Ga))$.

\bprp\label{soneggg}
$\ufin(\Op,\Ga)\cap\sone(\Ga,\Op)=\mbox{Lindel\"of }\sone(\Ga,\gimel(\Ga))$.
\eprp
\bpf
$(\Impl)$ $\ufin(\Op,\Ga)$ implies that every countable open cover is in $\gimel(\Ga)$ \cite{coc7}.

$(\Leftarrow)$ It suffices to prove that $\sone(\Ga,\gimel(\Ga))$ implies $\ufin(\Ga,\Ga)$.
Assume that $X$ satisfies  $\sone(\Ga,\gimel(\Ga))$, and let $\cU_n\in\Ga$ for all $n$.
We may assume that the covers $\cU_n$ get finer with $n$.
Apply $\sone(\Ga,\gimel(\Ga))$ to obtain $U_n\in\cU_n$, $n\in\N$, and
an increasing $h\in\NN$ such that, for each $x\in X$,
$$x\in \Un_{k=h(n)}^{h(n+1)-1}U_k$$
for all but finitely many $n$.
Since the covers $\cU_n$ get finer with $n$, there is for each $n$
a finite set $\cF_{h(n)}\sub \cU_{h(n)}$ such that
$$\Un_{k=h(n)}^{h(n+1)-1}U_k\sub\Un\cF_{h(n)}.$$
For $n$ not in the image of $h$ chose $\cF_n=\emptyset$.
\epf

\blem\label{addggimel}
Let a space $X$ be a union of less than $\cov(\cM)$ many $\sone(\Ga,\gimel(\Ga))$ spaces.
Then $X$ satisfies $\sone(\Ga,\Op)$.
\elem
\bpf
The proof is similar to that of Theorem 2.2 in \cite{AddGN}. We provide it, with the necessary changes, for completeness.

Let $\kappa<\cov(\cM)$.
Assume that, for each $\alpha<\kappa$, $X_\alpha$ satisfies $\sone(\Ga,\gimel(\Ga))$, and
$X=\Un_{\alpha<\kappa}X_\alpha$.
Let $\cU_n\in\Ga(X)$ for all $n$. We may assume that each $\cU_n$ is countable, and enumerate $\cU_n=\set{U^n_m}{m\in\N}$.
For each $\alpha$, as $X_\alpha$ satisfies $\sone(\Ga,\gimel(\Ga))$,
there are $f_\alpha\in\NN$ and an increasing $h_\alpha\in\NN$ such that,
for each $x\in X_\alpha$,
$$x\in \Un_{k=h_\alpha(n)}^{h_\alpha(n+1)-1}U^k_{f_\alpha(k)}$$
for all but finitely many $n$.

Since $\kappa<\cov(\cM)\le\fd$ \cite{BlassHBK}, there is an increasing $h\in\NN$
such that, for each $\alpha<\kappa$, the set
$$I_\alpha = \set{n}{[h_\alpha(n),h_\alpha(n+1))\sub [h(n),h(n+1))}$$
is infinite \cite{BlassHBK}.
For each $\alpha<\kappa$, define
$$g_\alpha\in \prod_{n\in I_\alpha} \N^{[h(n),h(n+1))}$$
by
$g_\alpha(n)=f_\alpha\rest[h(n),h(n+1))$ for all $n\in I_\alpha$.
As $\kappa<\cov(\cM)$, by Lemma 2.4.2(3) in \cite{BarJu},
there is $g\in\prod_n \N^{[h(n),h(n+1))}$ guessing all functions $g_\alpha$,
that is, for each $\alpha<\kappa$, $g(n)=g_\alpha(n)$ for infinitely many $n\in I_\alpha$ \cite{BlassHBK}.
Define $f\in\NN$ by $f(k)=g(n)(k)$, where $n$ is the one with $k\in[h(n),h(n+1))$.
Then $\set{U^n_{f(n)}}{n\in\N}\in\Op(X)$.

Indeed, let $x\in X$. Pick $\alpha<\kappa$ with $x\in X_\alpha$.
Pick $m$ such that, for all $n>m$,
$x\in \Un_{k=h_\alpha(n)}^{h_\alpha(n+1)-1}U^k_{f_\alpha(k)}$.
Pick $n\in I_\alpha$ such that $n>m$ and $g(n)=g_\alpha(n)$. Then
$$x\in \Un_{k=h_\alpha(n)}^{h_\alpha(n+1)-1}U^k_{f_\alpha(k)}\sub
\Un_{k=h(n)}^{h(n+1)-1}U^k_{f_\alpha(k)}=\Un_{k=h(n)}^{h(n+1)-1}U^k_{f(k)}.\qedhere$$
\epf

\bcor
$\add(\cM)\le\add(\sone(\Ga,\gimel(\Ga)))\le\fb$.
\ecor
\bpf
The second inequality follows from $\non(\ufin(\Op,\Ga))=\fb$.

First inequality: $\add(\cM)=\min\{\fb,\cov(\cM)\}$, and $\add(\ufin(\Op,\Ga))=\fb$.
Apply Lemma \ref{addggimel}.
\epf

We obtain the following generalization of Theorem \ref{concgggo}.

\bthm\label{concgimel}
Let $C$ be a $\cov(\cM)$-concentrated space. For each space $Y$ satisfying $\ufin(\Op,\Ga)$ and $\sone(\Ga,\Op)$,
$X\x Y$ satisfies $\sone(\Ga,\Op)$. \qed
\ethm
\bpf
Similar to the proof of Theorem \ref{concgggo}, using Lemma \ref{addggimel}.
\epf

Similarly, we have the following.

\bthm\label{conc4}
The product of each member of $\scrC_{\add(\cN)}(\cov(\cM))$ with every Lindel\"of $\sone(\Ga,\allowbreak\gimel(\Ga))$ space
satisfies $\sone(\Ga,\Op)$.\qed
\ethm

By Lemma \ref{addggimel}, $\cov(\cM)\le\non(\productive{\sone(\Ga,\gimel(\Ga))}{\sone(\Ga,\Op)})$.
Theorem \ref{conc4} holds, more generally, for $\scrC_{\add(\sone(\Ga,\Op))}(\non(\productive{\mbox{Lindel\"of }\sone(\Ga,\gimel(\Ga))}{\sone(\Ga,\Op)}))$.
Under mild hypotheses on a family $\scrA$ of covers,
the results proved here apply to $\sone(\scrA,\Op)$ for all $\scrA$.
The hypotheses on $\scrA$, which can be extracted from the proofs, are satisfied by all major types
of covers in the context of selection principles.

\bprb
Is $\add(\sone(\Ga,\gimel(\Ga)))=\fb$?
\eprb

\section{Concentrated sets and coherence of filters}

For $a\in\roth$ and an increasing $h\in\NN$, define
$$a/h = \set{n}{a\cap [h(n),h(n+1))\neq\emptyset}.$$
For $S\sub\roth$, define $S/h = \set{a/h}{a\in S}$.
%For an increasing $h\in\NN$ and a filter $\cF\sub\roth$,
%$$\cF/h= \{a/h : a\in\cF\} = \left\{ a : \Union_{n\in a}\intvl{h}{n}\in\cF\right\}.$$
%(The first equality is the definition; the second an easy fact.)

\subsection{Assuming NCF}
\emph{NCF (near coherence of filters)} is the assertion that, for each pair of nonprincipal ultrafilters
$\cU$ and $\cV$, there is an increasing $h\in\NN$ such that $\cU/h=\cV/h$.
The basic facts about NCF used here are available, e.g., in \cite{NCF1, BlassMil99}.

Henceforth, we use the convenient notation
$$U^n_{\le g(n)}:=\Un_{m=1}^{g(n)}U^n_m.$$

\bthm[NCF]\label{thm:Kd}
For a space $X$, the following assertions are equivalent:
\be
\item $X$ satisfies $\ufin(\Op,\Om)$.
\item Whenever $X\sub G\sub K$, with $K$ compact and $G$ G$_\delta$ in $K$,
there are $\kappa<\fd$ and compact sets $K_\alpha\sub K$, $\alpha<\kappa$,
such that $X\sub \Un_{\alpha<\kappa}K_\alpha\sub G$.
\ee
Moreover, the implication $(2)\Impl (1)$ holds in ZFC, and in the implication $(1)\Impl (2)$, we may take $\kappa=\fu$.
\ethm
\bpf
$(1)\Impl (2)$: Since $X$ is Lindel\"of, for each open set $U$ containing $X$, there are open sets $U_m$, $m\in\N$,
such that $X\sub\Un_m U_m\sub\Un_m\cl{U_m}\sub U$.
Let $G=\bigcap_n U_n$ with each $U_n$ open in $K$.
For each $n$, let $U^n_m$, $m\in\N$, be such that
$$X\sub\Un_{m\in\N} U^n_m\sub\Un_{m\in\N}\cl{U^n_m}\sub U_n,$$
and such that the covers $\set{U^n_m}{m\in\N}$ of $X$ get finer with $n$.
We may assume that none of these covers contains a finite subcover of $X$.

Apply $\ufin(\Op,\Om)$ to obtain an increasing $f\in\NN$ such that
$\set{U^n_{\le f(n)}}{n\in\N}\in\Om(X)$.
For each $x\in X$, let
$$f_x(n)=\min\set{m\ge f_x(n-1)}{x\in U^n_m}$$
for all $n$.
The family of all sets $\set{n}{f_x(n)\le f(n)}$, $x\in X$, is centered.
Extend it to a nonprincipal ultrafilter $\cU$ on $\N$.
Then $\set{f_x}{x\in X}$ is $\le_\cU$-bounded.

Let $\cV$ be an ultrafilter with base of size $\fu$.
By NCF, there is an increasing $h\in\NN$ such that $\cU/h=\cV/h$.
We claim that $\set{f_x}{x\in X}$ is $\le_\cV$-bounded.
Indeed, $f$ is an increasing $\le_\cU$-bound for  $\set{f_x}{x\in X}$.
Define $g(n)=f(h(n+1))$ for all $n$. For each $n\in \set{k}{f_x(k)\le f(k)}/h$, fix
$k\in \set{k}{f_x(k)\le f(k)}\cap [h(n),h(n+1))$, then
$$f_x(n)\le f_x(h(n))\le f_x(k)\le f(k)\le f(h(n+1))=g(n).$$
Let $\set{A_\alpha}{\alpha<\fu}$ be a base for $\cV$.
For each $\alpha<\fu$, let
$$K_\alpha=\bigcap_{n\in A_\alpha}\cl{U^n_{\le g(n)}}.$$
Then each $K_\alpha$ is compact, and $X\sub\Un_{\alpha<\fu}K_\alpha$.
NCF implies that $\fu<\fd$.

$(2)\Impl(1):$
For each $n$, let $\set{U^n_m}{m\in\N}$ be an open cover of $X$. Let $K$ be a compact space containing $X$.
We may assume that each set $U^n_m$ is open in $K$.
Let
$$G=\bigcap_{n\in\N}\bigcup_{m\in\N}U^n_m.$$
Then $G$ is G$_\delta$ in $K$, and $X\sub G$.
Let $\kappa<\fd$ and $K_\alpha\sub K$, $\alpha<\kappa$, be compact sets with
$X\sub \Un_{\alpha<\kappa}K_\alpha\sub G$.
For each $\alpha<\kappa$, let $f_\alpha\in\NN$ be such that
$$K_\alpha\sub \bigcap_n U^n_{\le f_\alpha(n)}.$$
Let $g\in\NN$ be a witness that $\set{f_\alpha}{\alpha<\kappa}$ is not finitely dominating.
Then $\set{U^n_{\le g(n)}}{n\in\N}$ is in $\Omega(X)$.
\epf

In \cite{huremen2} it was proved that, if NCF holds, then
$\fb,\fg\le \add(\ufin(\Op,\Om))$. We obtain an optimal version of this result.

\bcor[NCF]\label{cor:addSch}
$\add(\ufin(\Op,\Om))=\fd$.
\ecor
\bpf
By Theorem \ref{thm:Kd}.
\epf

\bthm[NCF]\label{thm:NCF}
Let $C$ be a $\fd$-concentrated space. For each $\ufin(\Op,\Om)$ space $Y$, $X\x Y$ satisfies $\sfin(\Op,\Op)$.
\ethm
\bpf
We use the projection method.

Assume that there is a countable $D\sub C$ with $|C\sm U|<\fd$ for all open $U\spst D$.
Let $Y$ be a  $\ufin(\Op,\Om)$ space.
Let $K$ be a compact space containing $C$ as a subspace.
Let $\cU_n$, $n\in\N$, be covers of $C\x Y$ by sets open in $K\x Y$.

As $D\x Y$ satisfies $\sfin(\Op,\Op)$, there are finite sets $\cF_n\sub\cU_n$, $n\in\N$, such that $D\x Y\sub U:=\Un_n\Un\cF_n$.

The projection $S$ of $(K\x Y)\sm U$ on the first coordinate, satisfies $\ufin(\Op,\Om)$.
By Theorem \ref{thm:Kd}, there are compact sets $K_\alpha$, $\alpha<\fu$, such that
$$S\sub \Un_{\alpha<\fu}K_\alpha\sub K\sm D.$$
As $C$ is $\fd$-concentrated on $D$, $|C\cap K_\alpha|<\fd$ for all $\alpha<\fu$.
By NCF, $\fu<\fd$ and $\fd$ regular. It follows that
$$|C\cap \Un_{\alpha<\fu}K_\alpha|<\fd.$$
By Corollary \ref{cor:addSch},
$(C\cap S)\x Y\sub (C\cap \Un_{\alpha<\fu}K_\alpha)\x Y$
and the latter space satisfies $\ufin(\Op,\Om)$. In particular, there are finite $\cF_n'\sub\cU_n$,
$n\in\N$, such that $(C\cap S)\x Y\sub \Un_n \Un\cF_n'$.
Then
$$C\x Y\sub\Un_{n\in\N}\Un(\cF_n\cup \cF_n').\qedhere$$
\epf

\subsection{Assuming $\fu<\fg$}
The axiom $\fu<\fg$ is stronger than NCF \cite{MdSh:894}.

\bthm[$\fu<\fg$]\label{thm:Ku}
For a space $X$, the following assertions are equivalent:
\be
\item $X$ satisfies $\sfin(\Op,\Op)$.
\item Whenever $X\sub G\sub K$, with $K$ compact and $G$ G$_\delta$ in $K$,
there are compact sets $K_\alpha\sub K$, $\alpha<\fu$,
such that $X\sub \Un_{\alpha<\fu}K_\alpha\sub G$.
\ee
\ethm
\bpf
$(2)\Impl (1)$: As $\fu<\fg\le\fd$, Theorem \ref{thm:Kd} implies that $X$ satisfies $\ufin(\Op,\Om)$.

$(1)\Impl (2)$: $\fu<\fg$ implies that $\ufin(\Op,\Om)=\sfin(\Op,\Op)$ \cite{SF1} (cf.\ \cite{SFT}).
Apply Theorem \ref{thm:Kd}.
\epf

\bcor[$\fu<\fg$]\label{cor:addMen}
$\add(\sfin(\Op,\Op))=\fd$.
\ecor
\bpf
By Corollary \ref{cor:addSch}, using that $\fu<\fg$ implies that $\ufin(\Op,\Om)=\sfin(\Op,\Op)$.
\epf

\bdfn\label{dfn:Kd}
Let $\scrK_\fd$ be the smallest (with respect to inclusion) class of topological spaces with the following properties:
\be
\itm Every singleton space is in $\scrK_\fd$.
\itm $\scrK_\fd$ is closed under unions of less than $\fd$ elements.
\itm If there is a countable $D\sub C$ with $C\sm U\in\scrK_\fd$ for all open $U\spst D$, then $C\in\scrK_\fd$.
\ee
\edfn

Notice that every $\fd$-concentrated space is in $\scrK_\fd$.

\bthm[$\fu<\fg$]\label{ugthm}
Every member of $\scrK_\fd$ is productively $\sfin(\Op,\Op)$ and productively $\sfin(\Om,\Om)$.
\ethm
\bpf
It suffices to prove the first assertion.
We use the projection method, and argue by induction on the structure of $\scrK_\fd$, as defined in Definition
\ref{dfn:Kd}. Case (1) in this definition is trivial, and Case (2) follows from Corollary \ref{cor:addMen}.
We treat Case (3).

Assume that there is a countable $D\sub C$ with $C\sm U$ productively $\sfin(\Op,\Op)$  for all open $U\spst D$.
Let $Y$ be a  $\sfin(\Op,\Op)$ space.
Let $K$ be a compact space containing $C$ as a subspace.
Let $\cU_n$, $n\in\N$, be covers of $C\x Y$ by sets open in $K\x Y$.

As $D\x Y$ satisfies $\sfin(\Op,\Op)$, there are finite sets $\cF_n\sub\cU_n$, $n\in\N$, such that $D\x Y\sub U:=\Un_n\Un\cF_n$.

The projection $M$ of $(K\x Y)\sm U$ on the first coordinate, satisfies $\sfin(\Op,\Op)$.
By Theorem \ref{thm:Ku}, there are compact sets $K_\alpha$, $\alpha<\fu$, such that
$$M\sub \Un_{\alpha<\fu}K_\alpha\sub K\sm D.$$
Let $\alpha<\fu$. As $K_\alpha\cap D=\emptyset$, we have by the induction hypothesis that
$(C\cap K_\alpha)\x Y$ satisfies $\sfin(\Op,\Op)$.
By Corollary \ref{cor:addMen},
$$\Un_{\alpha<\fu}(C\cap K_\alpha)\x Y\sub (C\cap M)\x Y$$
satisfies $\sfin(\Op,\Op)$.
Take finite $\cF_n'\sub\cU_n$, $n\in\N$, such that $(C\sm G)\x Y\sub \Un_n \Un\cF_n'$.
Then
$$C\x Y\sub\Un_{n\in\N}\Un(\cF_n\cup \cF_n').\qedhere$$
\epf

A notorious open problem asks whether, consistently, $\sfin(\Op,\Op)$ is closed under finite products.
By Theorem \ref{ugthm}, a positive answer to the following problem would settle this problem in the affirmative.
The \emph{superperfect set model} is the model obtained by an $\aleph_2$ stage countable support iteration of superperfect trees
forcing over a model of GCH. In this model, $\fu<\fg$. The values of the combinatorial cardinal characteristics of the continuum in this model \cite{BlassHBK}
imply that there are no generalized (in any relevant sense) Luzin or Sierpi\'nski
sets there (see Section \ref{sec:LuzSier} for the definitions). Consequently, in the superperfect set model,
the only known spaces satisfying $\sfin(\Op,\Op)$ are those in $\scrK_\fd$.

\bprb\label{Miller}
Is $\scrK_\fd=\sfin(\Op,\Op)$ in the superperfect set model?
\eprb

We conclude this section with analogous results for Rothberger's property $\sone(\Op,\Op)$.
The hypothesis $\fu<\fg$ implies that every $\sone(\Op,\Op)$ space is $\ufin(\Op,\Ga)$ \cite{SF1}
(cf.\ \cite{SFT}), and therefore that $\sone(\Om,\gimel(\Ga))=\sone(\Op,\Op)$.

\bdfn\label{dfn:CcovM}
Let $\scrC_{\cov(\cM)}$ be the smallest (with respect to inclusion) class of topological spaces with the following properties:
\be
\itm Every singleton space is in $\scrC_{\cov(\cM)}$.
\itm $\scrC_{\cov(\cM)}$ is closed under unions of less than $\cov(\cM)$ elements.
\itm If there is a countable $D\sub C$ with $C\sm U\in\scrC_{\cov(\cM)}$ for all open $U\spst D$, then $C\in\scrC_{\cov(\cM)}$.
\ee
\edfn

Notice that every $\cov(\cM)$-concentrated space is in $\scrC_{\cov(\cM)}$.
Using the above methods, we obtain the following.

\bthm[$\fu<\fg$]\label{ugthm2}
Every member of $\scrC_{\cov(\cM)}$ is productively $\sone(\Op,\Op)$ and productively $\sone(\Om,\Om)$.
\ethm
\bpf
Assuming $\fu<\fg$, since $\sone(\Om,\gimel(\Ga))=\sone(\Op,\Op)$, we have by Theorem 2.3 of \cite{AddGN} that
\begin{eqnarray*}
\add(\sone(\Op,\Op)) & = & \add(\sone(\Om,\gimel(\Ga)))=\add(\cM)=\non(\sone(\Om,\gimel(\Ga)))=\\
& = & \non(\sone(\Op,\Op))=\cov(\cM).
\end{eqnarray*}
The rest follows from the projection method, as in the proof of Theorem \ref{ugthm}.
\epf

A notorious open problem, due to Bartoszy\'nski and Judah \cite{covM2}, asks whether, consistently, $\add(\sone(\Op,\Op))$ or $\add(\sone(\B,\B))$ may be greater than 
$\add(\cN)$.
We do not know whether the hypothesis in the following theorem is consistent
(it is provable that $\cov(\cM)\le\fu$, though), but once such a consistency result is established, we would obtain a solution of this problem.

\bthm[$\add(\cN)<\cov(\cM)\le\fu<\fg$]
$$\add(\cN)<\cf(\cov(\cM))=\cov(\cM)=\add(\sone(\Op,\Op))=\add(\sone(\B,\B)).$$
\ethm
\bpf
By the arguments in the proof of Theorem \ref{ugthm2} and the results used to prove it (all applying to $\sone(\B,\B)$ as well),
we have that
$$\add(\cM)=\cov(\cM)=\add(\sone(\Op,\Op))=\add(\sone(\B,\B)).$$
Since $\add(\cM)$ is regular, so is $\cov(\cM)$.
\epf

\bprb
Is it consistent that $\add(\cN)<\cov(\cM)\le\fu<\fg$?
\eprb

As in Problem \ref{Miller}, we do not know whether $\scrC_{\cov(\cM)}=\sone(\Op,\Op)$ in the superperfect set model, or
whether there is at all a model where Borel's Conjecture fails (i.e., there are uncountable $\sone(\Op,\Op)$ sets of reals)
and $\sone(\Op,\Op)$ is closed under finite products.

\section{Luzin and Sierpi\'nski sets}\label{sec:LuzSier}

Let $\kappa$ be an uncountable cardinal.
A set $L\sub\R$ is $\kappa$-\emph{Luzin} if its intersection with every meager subset of $\R$ has
cardinality less than $\kappa$. \emph{Luzin sets} are $\aleph_1$-Luzin subsets of $\R$.
Every $\kappa$-Luzin set $L$ is $\kappa$-\emph{concentrated} on a countable subset $D\sub L$ (indeed, on every
countable dense subset $D\sub L$).
A set $S\sub \R$ is $\kappa$-\emph{Sierpi\'nski} if  its intersection with every Lebesgue null subset of $\R$ has
cardinality less than $\kappa$. \emph{Sierpi\'nski sets} are $\aleph_1$-Sierpi\'nski subsets of $\R$.

The starting point of this section, that indeed also led to the earlier two sections,
is a surprising result of Babinkostova and Scheepers.\footnote{For
the following details, it is recommended to consult the Scheepers Diagram.}
Let $L$ be a Luzin set and $S$ be a Sierpi\'nski set. It is known (e.g., \cite{coc2, CBC}) that:
\be
\itm $L$ satisfies $\sone(\Op,\Op)$ (indeed, $\sone(\B,\B)$), but not $\ufin(\Op,\Ga)$.
\itm $S$ satisfies $\sone(\Ga,\Ga)$ (indeed, $\sone(\BG,\BG)$), but not $\sone(\Op,\Op)$.
\ee
Assuming \CH{}, there is a Luzin set $L$ that does not satisfy $\ufin(\Op,\Om)$ \cite{coc2}.
It follows that $L\x S$ does not satisfy any of $\ufin(\Op,\Om)$ or $\sone(\Op,\Op)$.

\brem
It follows, in particular, that one cannot improve Theorem \ref{concgggo} by proving, e.g.,
that every concentrated set of real numbers is productively $\sone(\Ga,\Ga)$.
\erem

\bthm[Babinkostova--Scheepers \cite{BabSch}]\label{BSThm}
For every Luzin set $L$ and Sierpi\'nski set $S$, $L\x S$ satisfies $\sfin(\Op,\Op)$.
\ethm

As mentioned in the earlier sections, Babinkostova and Scheepers prove in \cite{BabSch} that it
suffices to assume that $L$ is concentrated on a countable subset (or even less),
and that $S$ satisfies $\ufin(\Op,\Ga)$. Their full result is generalized
further in \cite{AddGN}. In the present section, we settle the question which
selective properties are provably satisfied by products of Luzin and Sierpi\'nski sets.
First, we use the results of the earlier section to settle the problem in the case of open covers.

\bthm
For every Luzin set $L$ and Sierpi\'nski set $S$, $L\x S$ satisfies $\sone(\Ga,\Op)$.
\ethm
\bpf
Recall that Luzin sets are concentrated on countable subsets, and Sierpi\'nski sets satisfy
$\sone(\Ga,\Ga)$. Apply Theorem \ref{concgimel} (or Theorem \ref{concgggo}).
\epf

With, apparently, no exceptions thus far, all results about selective covering properties of Luzin and Sierpi\'nski
sets, proved in the realm of open covers, were also provable for the corresponding Borel-covers variant.
Some examples are available in \cite{CBC}.
In light of this, the results in the remainder of this section are surprising.
They imply, in particular, that a product of a Luzin and a Sierpi\'nski
set may \emph{fail} to satisfy $\sfin(\B,\B)$ (Menger's property for Borel covers), and thus
any of the Borel-cover versions of the properties in the Scheepers diagram.

For convenience, in the remainder of this section we work in the Cantor space $\Cantor$ instead of $\R$.
The results can be transformed into $\R$ using the canonical map
\begin{eqnarray*}
\Cantor & \to & [0,1]\\
f & \mapsto & \sum_n \frac{f(n)}{2^n}.
\end{eqnarray*}

\bdfn
Define a reflexive binary relation $R$ on $\Cantor$ by setting $xRy$ if
$$\exists^\infty n,\ x\upharpoonright [n,2n)=y\upharpoonright[n,2n).$$
For $y\in \Cantor$, let
$$[y]_R=\set{x\in \Cantor}{xRy}.$$
For $\vec x=\seq{x_n}{n\in\N}\in\cube$ and $y\in \Cantor$, define
$$\Match(\vec x,y)=\chi_{\set{n}{x_nRy}}.$$
\edfn

\blem\label{basic1}
For $x,y\in \Cantor$:
\be
\itm $[y]_R$ is a Lebesgue null, $G_\delta$ dense subset of $\Cantor$.
\itm If $x=^*y$ (equal mod finite), then $[x]_R=[y]_R$.
\itm $\Match\colon\cube\x \Cantor\to \Cantor$ is a Borel map.
\qed
\ee
\elem

\bthm
For all comeager $X\sub\cube$ and nonnull $Y\sub \Cantor$, $\Match[X\x Y]=\Cantor$.
\ethm
\bpf
\blem\label{1}
Let $Y\sub \Cantor$ be nonnull.
For each $\vec x\in\cube$, there is $y\in Y$ such that $\vec x\in({[y]_R}\comp)^\N$.
\elem
\bpf
By Lemma \ref{basic1}, $[x_n]_R$ is null for all $n$, and consequently so is $\Un_n[x_n]_R$.
Pick
$$y\in Y\sm \Un_n[x_n]_R.$$
As $R$ is symmetric, $\vec x\in({[y]_R}\comp)^\N$.
\epf
\blem\label{2}
Let $Y\sub \Cantor$ be nonnull.
For each $I\sub\N$, $\Un_{y\in Y}({[y]_R}\comp)^{I\comp}\x ([y]_R)^I$ is nonmeager in $\cube$.
\elem
\bpf
Let $\vec x\in (\Cantor)^{I\comp}$.
By Lemma \ref{1}, there is $y\in Y$ such that $\vec x\in({[y]_R}\comp)^{I\comp}$.
By Lemma \ref{basic1}, $([y]_R)^I$ is comeager in $(\Cantor)^I$.
Since
$$\{\vec x\}\x ([y]_R)^I\sub \Un_{y\in Y}({[y]_R}\comp)^{I\comp}\x ([y]_R)^I,$$
it follows that all vertical sections of $\Un_{y\in Y}({[y]_R}\comp)^{I\comp}\x ([y]_R)^I$ are comeager,
in particular nonmeager.
By \cite[Lemma 8.42]{Ke}, our set is nonmeager in $(\Cantor)^{I\comp}\x (\Cantor)^I$.
\epf
Let $z=\chi_I\in \Cantor$.
By Lemma \ref{2}, there is
$$\vec x\in X\cap\left(\Un_{y\in Y}({[y]_R}\comp)^{I\comp}\x ([y]_R)^I\right).$$
Then $\Match(\vec x,y)=\chi_I=z$.
\epf

\bcor\label{BorelCantor}
There is a Borel map $f\colon \Cantor\to \Cantor$ such that, for all comeager $X\sub \Cantor$ and nonnull $Y\sub \Cantor$,
$f[X\x Y]=\Cantor$.
\ecor
\bpf
The canonical bijection $\Cantor\to\cube$ is Borel, and preserves meager and null sets in both directions.
\epf

Sierpi\'nski sets are special kinds of nonnull sets.
In the Sacks model, there are Luzin and Sierpi\'nski sets, but they are all of cardinality $\aleph_1$, whereas the continuum
is $\aleph_2$. Thus, consistently, there are no Luzin and Sierpi\'nski sets whose product can be mapped onto $\Cantor$.
However, we have the following.

\bcor[CH]
For each nonnull set $Y\sub \Cantor$, there is a Luzin set $L\sub \Cantor$ such that $\Cantor$ is a Borel image of $L\x Y$.
\ecor
\bpf
Let $f$ be the function defined in Corollary \ref{BorelCantor}.
Enumerate $\Cantor=\set{r_\alpha}{\alpha<\aleph_1}$.
Let $\set{M_\alpha}{\alpha<\aleph_1}$ be a cofinal family of meager subsets of $\Cantor$.
For each $\alpha<\aleph_1$,
$$f\left[ \left(\Cantor\sm \Un_{\beta<\alpha}M_\alpha\right)\x Y \right]=\Cantor.$$
Pick $(x_\alpha,y_\alpha)\in (\Cantor\sm \Un_{\beta<\alpha}M_\alpha)\x Y$ such that
$f(x_\alpha,y_\alpha)=r_\alpha$.
Finally, let $L=\set{x_\alpha}{\alpha<\aleph_1}$.
\epf

$\add(\cM)$-Luzin sets satisfy $\sone(\B,\B)$ (e.g., \cite{CBC}).

\bcor[$\add(\cM)=\cof(\cM)$]\label{happy}
For each nonnull set $Y\sub \Cantor$, there is an $\add(\cM)$-Luzin set $L\sub \Cantor$ such that $L\x Y$ does not satisfy $\sfin(\B,\B)$.
\ecor
\bpf
Let $f$ be the function defined in Corollary \ref{BorelCantor}.
Let $\kappa=\add(\cM)=\cof(\cM)$.
As $\add(\cM)\le\fd\le\cof(\cM)$ (in ZFC), there is a dominating set $\set{d_\alpha}{\alpha<\kappa}\sub\NN$.
Identify $\Cantor$ with $\NN$ via a Borel bijection.

Let $\set{M_\alpha}{\alpha<\kappa}$ be a cofinal family of meager subsets of $\Cantor$.
For each $\alpha<\aleph_1$,
$$f\left[ \left(\Cantor\sm \Un_{\beta<\alpha}M_\alpha\right)\x Y \right]=\NN.$$
Pick $(x_\alpha,y_\alpha)\in (\Cantor\sm \Un_{\beta<\alpha}M_\alpha)\x Y$ such that
$f(x_\alpha,y_\alpha)=d_\alpha$.
Finally, let $L=\set{x_\alpha}{\alpha<\aleph_1}$.

As the Borel image $f[L\x Y]$ contains $D$, it is dominating. Thus, $L\x Y$ does not satisfy $\sfin(\B,\B)$.
\epf

\bthm[$\add(\cN)=\cof(\cN)$]
There are an $\add(\cN)$-Luzin set $L\sub \Cantor$ and an $\add(\cN)$-Sirepi\'nski set $S\sub \Cantor$ such that:
\be
\itm All finite powers of $L$ satisfy $\sone(\BO,\BO)$;
\itm All finite powers of $S$ satisfy $\sone(\BG,\BG)$ and $\sfin(\BO,\BO)$; but
\itm $L\x S$ does not satisfy $\sfin(\B,\B)$.
\ee
Moreover, $L$ does not satisfy $\ufin(\Op,\Ga)$ and $S$ does not satisfy $\sone(\Op,\Op)$.
\ethm
\bpf
As $\add(\cN)=\cof(\cN)$, there are a $\add(\cN)$-Sierpi\'nski set $S$ as in (2) \cite[Corollary 25]{o-bdd},
and by a dual argument,
an $\add(\cN)$-Luzin set $L$ such that all finite powers of $L$ satisfy $\sone(\B,\B)$.
(Here, we use Carlson's Theorem, that the union of less than $\add(\cN)$ elements of $\sone(\B,\B)$
is in $\sone(\B,\B)$ \cite{AddQuad}.)
It is pointed out in \cite[Theorem 18]{CBC} that
if all finite powers of $X$ have property $\sone(\B,\B)$, then $X$ has property $\sone(\BO,\BO)$.
This implies (1).

During the construction of $L$, one can also accommodate the restrictions provided in the proof of Corollary \ref{happy}, to make sure that $f[L\x S]$ contains a (Borel preimage in $\Cantor$ of a)
dominating subset of $\NN$.
This gives (3).

The last assertion in the theorem is due to Sierpi\'nski, cf.\ \cite{coc2}.
\epf

\section{Scales and $\fb$-scales}\label{sec:scales}

In the earlier sections, we have discussed Luzin sets as special examples of concentrated sets.
Another standard method for constructing concentrated sets, initiated by Rothberger, is that of using \emph{scales}.
These constructions require in general milder hypotheses than those used for the construction of Luzin and Sierpi\'nski sets,
and in many cases can be carried out outright in ZFC.

For our purposes, it is convenient to identify the Cantor space $\Cantor$ with $\PN=\roth\cup\Fin$,
and use the induced topology.
For $a\in\roth$ and $n\in\N$, $a(n)$ denotes the $n$-th element in the increasing enumeration of $a$.
For $a,b\in\roth$, $a\le^* b$ means that $a(n)\le b(n)$ for all but finitely many $n$.
%$\fb,\fd$ denote, respectively, the minimal cardinalities of an unbounded and a cofinal (or dominating)
%subset of $\roth$, with respect to $\le^*$.
A \emph{scale} is a cofinal (dominating) set $S=\set{s_\alpha}{\alpha<\kappa}$ in $(\roth,\le^*)$ such that
$s_\alpha\le^* s_\beta$ for $\alpha<\beta$.
Scales exist if and only if $\fb=\fd$, and in this case, their cardinality is $\fb$.
If we generalize ``cofinal'' to ``unbounded'', we obtain the definition of \emph{$\fb$-scale},
an object constructible within ZFC.
For each $\fb$-scale $B$, $B\cup\Fin$ is $\fb$-concentrated on its countable subset $\Fin$.
%We will introduce an additional important variant of scales later on.

For brevity, the union of a scale with $\Fin$, viewed as a subset of the Cantor space $\PN$,
will be called \emph{scale set}. \emph{$\fb$-scale sets} are defined similarly.

\bcor
Every product of a $\fb$-scale set and a $\cov(\cM)$-concentrated space satisfies $\sone(\Ga,\Op)$.
\ecor
\bpf
Every $\fb$-scale set satisfies $\ufin(\Op,\Ga)$ and, being $\fd$-concentrated, $\sone(\Ga,\Op)$ as well \cite{ideals}.
Apply Proposition \ref{soneggg} and Theorem  \ref{concgimel}.
\epf

\bthm\label{scoo}
Let $S\sub\roth$ be a scale. The scale set $S\cup\Fin$ is productively $\sfin(\Op,\Op)$ and productively $\sfin(\Om,\Om)$.
\ethm
\bpf
Since $\sfin(\Om,\Om)=\sfin(\Op,\Op)^\uparrow$, it suffices to prove the first assertion.
We use the projection method (cf.\ Theorem \ref{concgggo}).

Let $Y$ be a space satisfying $\sfin(\Op,\Op)$.
Let $\cU_n$, $n\in\N$, be covers of $(S\cup\Fin)\x Y$ by sets open in $\PN\x Y$.
As $\Fin$ is countable, $\Fin\x Y$ satisfies $\sfin(\Op,\Op)$.
Pick finite sets $\cF_n\sub\cU_n$, $n\in\N$, such that $\Fin\x Y\sub U:=\Un_n\Un\cF_n$.

Since $Y$ satisfies $\sfin(\Op,\Op)$ and $\PN$ is compact, the projection $M$ of $(\PN\x Y)\sm U$
on the first coordinate satisfies $\sfin(\Op,\Op)$. $M\sub\roth$, and satisfying $\sfin(\Op,\Op)$, it is not dominating.
Thus, $|M\cap S|<\fd=\fb$. It follows that $(M\cap S)\x Y$ satisfies $\sfin(\Op,\Op)$.
Pick finite sets $\cF_n'\sub\cU_n$, $n\in\N$, such that $(M\cap S)\x Y\sub \Un_n\Un\cF_n'$.
Then
$$X\x Y\sub \Un_n\Un(\cF_n\cup\cF_n').\qed$$
\epf

Theorem \ref{scoo} is the last one in this paper proved by the projection method.
In order to establish additional productive properties of scale sets, we use the following method.

\blem[Productive Two Worlds Lemma]\label{twl}
Let $Y$ be a space, and for each $n$, let
$$\set{U^n_m}{m\in\N}\in\Om(\Fin\x Y)$$
with each $U^n_m$ clopen in $\PN\x Y$.
There is a continuous map $\Psi\colon Y\to\NN$
such that, for each $n$,
$$(x,y)\in U^n_{\le\Psi(y)(n)}$$
for all $x\in\roth$ such that $|x|<n$ or $\Psi(y)(n)\le x(n)$.
\elem
\bpf
Let $y\in Y$.
Fix $n$. Let $m_n(1)=a_n(1)=1$.
By induction on $k$, let $m_n(k+1)$ be minimal with
$$P(\{1,\dots,a_n(k)-1\})\x\{y\}\sub U^n_{m_n(k+1)},$$
and let $a_n(k+1)$ be minimal such that
$$(x,y)\in U^n_{m_n(k+1)}$$
for all $x\in\PN$ with $x\cap \{a_n(k),\dots,a_n(k+1)-1\}=\emptyset$.

Define
$$\Psi(y)(n)=\max\{a_n(n+1),m_n(n+1)\}$$
for all $n$. $\Psi$ is continuous.
Fix $n$, and let $x\in\roth$ with $\Psi(y)(n)\le x(n)$.
As $a_n(n+1)\le \Psi(y)(n)\le x(n)$, there is $k<n+1$ with $x\cap \{a_n(k),\dots,a_n(k+1)-1\}=\emptyset$.
Then $(x,y)\in U^n_{m_n(k+1)}$. As
$$m_n(k+1)\le m_n(n+1)\le \Psi(y)(n),$$
$U^n_{m_n(k+1)}\sub U^n_{\le\Psi(y)(n)}$.
\epf

\blem\label{0dimlemma}
Let $X$ be a space and $\scrA\in \{\Ga,\Om,\Op\}$.
If there is a space $Y$ such that $X\x Y$ is  Lindel\"of but not $\ufin(\Op,\scrA)$,
then there such a subspace $Y$ of the Cantor space.
\elem
\bpf
Let $\cU_n$, $n\in\N$, be  open covers of $X\x Y$ witnessing the failure of $\ufin(\Op,\scrA)$.
As $X\x Y$ is Lindel\"of, we may assume that each $\cU_n$ has the form $\set{U^n_m \x V^n_m}{m\in\N}$.
Define set-valued maps $\Phi,\Psi$ from $X,Y$, respectively, into the Cantor space $P(\N\x\N)$ by
\begin{eqnarray*}
A\in \Phi(x) &\Longleftrightarrow& \set{(n,m)}{x\in U^n_m}\sub A\\
A\in \Psi(y) &\Longleftrightarrow& \set{(n,m)}{y\in V^n_m}\sub A
\end{eqnarray*}
for all $x\in X,y\in Y$.
By \cite[Lemma 2]{SF1}, these maps are compact-valued, upper semicontinuous.
Thus,
$$\Psi[Y]:=\Un_{y\in Y}\Psi(y)\sub P(\N\x\N)$$
satisfies $\ufin(\Op,\scrA)$, and $X\x\Psi[Y]$ is Linde\"of, being a compact-valued, upper semicontinuous
image of the Lindel\"of space $X\x Y$.

We claim that $X\x\Psi[Y]$ does not satisfy $\ufin(\Op,\scrA)$.
Assume otherwise. Then $\Phi[X]\x\Psi[Y]$ satisfies $\ufin(\Op,\scrA)$, being a
compact-valued, upper semicontinuous image of $X\x\Psi[Y]$.
Define
\begin{eqnarray*}
\Xi\colon  P(\N\x\N)\x  P(\N\x\N) & \to & P(\N\x\N)\\
(A,B) & \mapsto & A\cap B.
\end{eqnarray*}
Then $\Xi$ is continuous.
For each $n$, let $\cW_n:=\set{W^n_m}{m\in\N}$,
where
$$W^n_m=\set{A\sub \N\x\N}{(n,m)\in  A}.$$
For each $n$, $\cW_n$ is an open cover of $\Xi(\Phi[X]\x\Psi[Y])$.
But $\set{W^n_{\le f(n)}}{n\in\N}$ is not in $\scrA(\Xi(\Phi[X]\x\Psi[Y]))$
for any $f\in\NN$, for otherwise, $\set{\Un_{m\leq f(n)}U^n_m\x V^n_m}{n\in\N}$ would be in $\scrA(X\x Y)$.
\epf

An open cover $\cU\in\Omega^\mathrm{gp}(X)$ if there are $h\in\NN$ and an enumeration $\cU=\set{U_n}{n\in\N}$ such that, for each finite $F\sub X$ and each $n$,
$F\sub U_k$ for some $h(n)\le k\le h(n+1)$. $\ufin(\Op,\Ga)^\uparrow=\sfin(\Om,\Omega^\mathrm{gp})$ \cite{coc7}.

Bartoszy\'nski and Shelah \cite{BaSh01} proved that every $\fb$-scale set satisfies $\ufin(\Op,\Ga)$. Then,
Bartoszy\'nski and Tsaban \cite{ideals} proved that all finite powers of a $\fb$-scale set satisfy $\ufin(\Op,\Ga)$.
Later, Tsaban and Zdomskyy \cite{sfh} proved that all finite products of $\fb$-scale sets satisfy $\ufin(\Op,\Ga)$.
The following theorem is much stronger.

\bthm\label{bscalethm}
Every $\fb$-scale set (in particular, every scale set) is:
\be
\itm Productively hereditarily Lindel\"of $\ufin(\Op,\Ga)$; and
\itm Productively hereditarily Lindel\"of $\sfin(\Om,\Omega^\mathrm{gp})$.
\ee
\ethm
\bpf
It suffices to prove the first assertion.

Let $B=\set{b_\alpha}{\alpha<\fb}$ be a $\fb$-scale.
Let $Y$ be a hereditarily Lindel\"of space satisfying $\ufin(\Op,\Ga)$. Then $(B\cup\roth)\x Y$ is (hereditarily) Lindel\"of.
By Lemma \ref{0dimlemma}, we may assume that $Y\sub\Cantor$.

Let $\cU_n\in\Om((B\cup\Fin)\x Y)$ for all $n$. As $X\x Y$ is a subspace of the Cantor space,
we may assume that each element of each $\cU_n$ is clopen, and the Productive Two Worlds Lemma applies.
Let $\Psi$ be as in that Lemma. As $Y$ satisfies $\ufin(\Op,\Ga)$, $\Psi[Y]$ is bounded by some $g\in\NN$.
Take $\alpha<\fb$ such that the set $I=\set{n}{g(n)\le b_\alpha(n)}$ is infinite.
For each $\beta\ge \alpha$ and each $y\in Y$,
$$\Psi(y)(n)\le g(n)\le b_\alpha(n)\le b_\beta(n),$$
and therefore $(b_\beta,y)\in U^n_{\le g(n)}$, for all but finitely many $n\in I$.
We also have that  $(x,y)\in U^n_{\le g(n)}$, for all $x\in\Fin, y\in Y$ and all but finitely many $n$.

As $\add(\ufin(\Op,\Ga))=\fb$, $\set{b_\beta}{\beta<\alpha}\x Y$ satisfies $\ufin(\Op,\Ga)$,
and thus there is $h\in\NN$ such that $\set{U^n_{\le h(n)}}{n\in\N}\in\Ga(\set{b_\beta}{\beta<\alpha}\x Y)$.
For $n\in I$, let
$$\cF_n=\set{U^n_m}{m\le\max\{g(n),h(n)\}}.$$
For $n\notin I$, let $\cF_n=\emptyset$.
Then $\set{\Un\cF_n}{n\in\N}\in\Ga((B\cup\Fin)\x Y)$.
\epf

A set $D\sub\NN$ is \emph{finitely dominating} if its closure under pointwise maxima of finite subsets,
$\maxfin(D)$, is dominating. Let $\cov(\Dfin)$ be the minimal $\kappa$ such that $\NN$ (equivalently, a dominating
subset of $\NN$) can be decomposed into $\kappa$ many sets, none of which finitely dominating. Then
$$\max\{\fb,\fg\}\le\cov(\Dfin)\le\fd,$$
and strict inequalities are consistent \cite{MShT847}.

%Let $\Op_\mathrm{ctbl}(X)$ denote the family of all \emph{countable} open covers of $X$.

\blem\label{dfinlem}
Every space of cardinality less than $\cov(\Dfin)$ is productively $\ufin(\Op,\Om)$ for countable covers.
\elem
\bpf
Assume that $|X|<\cov(\Dfin)$ and $Y$ satisfies $\ufin(\Op,\Om)$ for countable covers.
Using the terminology of the forthcoming Section \ref{sec:chars} and Theorem \ref{thm:Schar},
let $\Psi\colon X\x Y\to\NN$ be upper continuous.
It suffices to prove that $\Psi[X\x Y]$ is not finitely dominating.
For each finite $F\sub X$, the map
\begin{eqnarray*}
\Psi_F\colon Y & \to & \NN\\
y & \mapsto & \max\set{\Psi(x,y)}{x\in F}
\end{eqnarray*}
is upper continuous. Thus, $\Psi_F[Y]$ is not finitely dominating.
Each finite subset of $X\x Y$ is contained in one of the form $F_1\x F_2$,
and
$$\max\Psi[F_1\x F_2]=\max\Psi_{F_1}[F_2].$$
Thus,
$$\Un_{F\in[X]^{<\infty}}\maxfin\Psi_{F}[Y]$$
is cofinal in $\maxfin\Psi[X\x Y]$.
This is a directed union (every finite sub-union is contained in a single member)
of less than $\cov(\Dfin)$ many sets that are not finitely dominating.
Thus, it is not finitely dominating.
\epf

\bthm
For each scale $S\sub\roth$, the scale set $S\cup\Fin$ is productively hereditarily Lindel\"of $\ufin(\Op,\Om)$.
\ethm
\bpf
Let $S=\set{s_\alpha}{\alpha<\fb}$ be a scale, and let $Y$ be a space satisfying $\ufin(\Op,\Om)$.
By Lemma \ref{0dimlemma}, we may assume that $Y$ is a subspace of the Cantor space.
Let
$$\cU_n=\set{U^n_m}{m\in\N}\in\Om((S\cup\Fin)\x Y)$$
for all $n$.
We may assume that every $U^n_m$, $n,m\in\N$, is clopen. Let $\Psi$ be as in Lemma \ref{twl}.
For $\alpha<\gamma<\fb$, let
$$X_{\alpha,\gamma}=\set{s_\beta}{\beta<\alpha\mbox{ or }\gamma\le\beta}\cup\Fin.$$
\blem\label{ag}
For each $\alpha<\fb$, there are $\gamma<\fb$ and $g\in\NN$ such that
$$\set{U^n_{\le g(n)}}{n\in\N}\in \Om\left(X_{\alpha,\gamma}\x Y\right).$$
\elem
\bpf
Fix $p\notin \set{s_\beta}{\beta<\alpha}$, and consider $D= \set{s_\beta}{\beta<\alpha}\cup\{p\}$ as a discrete
space.
By Lemma \ref{dfinlem}, $D\x Y$ satisfies $\ufin(\Op,\Om)$.
Define $\Phi\colon D\x Y\to\NN$ by
\begin{eqnarray*}
\Phi(p,y) & = & \Psi(y);\\
\Phi(s_\beta,y)(n) & = & \min\set{m}{(s_\beta,y)\in U^n_m}\quad (n\in\N)
\end{eqnarray*}
As $\Phi$ is continuous, $\Phi[D\x Y]$ is not finitely dominating. Let $g\in\NN$ be a witness for that.
Let $\gamma$ be such that $g\le^* s_\gamma$. We claim that $\gamma$ and $g$ are as required.
Let $F\sub X_{\alpha,\gamma}$ and $G\sub Y$ be finite sets. Decompose $F$ as
$$F=(F\cap X_{\alpha,\fb}\sm\Fin)\cup (F\cap\Fin)\cup (F\cap X_{0,\gamma}\sm\Fin).$$
The set
$$I:=\set{n}{\max\Phi[((F\cap X_{\alpha,\fb}\sm\Fin)\cup\{p\})\x G](n)\le g(n)}$$
is infinite.

For each $s_\beta\in F\cap X_{\alpha,\fb}\sm\Fin$ and each $y\in Y$,
$$\Phi(s_\beta,y)(n)\le g(n)$$
for all $n\in I$. Thus, $(F\cap X_{\alpha,\fb}\sm\Fin)\x G\sub U^n_{\le g(n)}$ for all $n\in I$.

For each $x\in\Fin$ and each $y\in G$,
$$\Psi(y)(n)=\Phi(p,y)(n)\le g(n)\mbox{ and }|x|<n$$
for all but finitely many $n\in I$. Thus, $(F\cap\Fin)\x G\sub U^n_{\le g(n)}$ for all but finitely many $n\in I$.

Finally, for each $\beta\ge\gamma$ and each $y\in G$,
$$\Psi(y)(n)=\Phi(p,y)(n)\le g(n)\le s_\beta(n)$$
for all but finitely many $n\in I$. Thus,  $(F\cap X_{0,\gamma}\sm\Fin)\x G\sub U^n_{\le g(n)}$ for all but finitely many $n\in I$.

It follows that $F\x G\sub U^n_{\le g(n)}$ for all but finitely many $n\in I$.
\epf
By Lemma \ref{ag} applied to $\alpha=0$, there are $\gamma_1<\fb$ and $g_1\in\NN$ such that
$$\set{U^n_{\le g_1(n)}}{n\in\N}\in \Om\left(X_{0,\gamma_1}\x Y\right).$$
By Lemma \ref{ag} applied to $\alpha=\gamma_1$, there are $\gamma_2<\fb$ and $g_2\in\NN$ such that
$$\set{U^n_{\le g_2(n)}}{n\in\N}\in \Om\left(X_{\gamma_1,\gamma_2}\x Y\right).$$
By Lemma \ref{ag} applied to $\alpha=\gamma_2$, there are $\gamma_3<\fb$ and $g_3\in\NN$ such that
$$\set{U^n_{\le g_3(n)}}{n\in\N}\in \Om\left(X_{\gamma_2,\gamma_3}\x Y\right).$$
Continue in, in this manner, to define $\gamma_n$ and $g_n$ for all $n$. Let $\gamma=\sup_n\gamma_n$, and
$g$ be a $\le^*$-bound of $\set{g_n}{n\in\N}$.
Then
$$\set{U^n_{\le g(n)}}{n\in\N}\in \Om\left((S\cup\Fin)\x Y\right).$$
Indeed, let $F\x G$ be a finite subset of $(S\cup\Fin)\x Y$. Let $H=\set{\beta<\fb}{s_\beta\in F}$.
As $H$ is finite, there is $k$ such that $H\cap\gamma\sub \gamma_k$. Then
$F\sub X_{\gamma_k,\gamma_{k+1}}$, and thus there are infinitely many $n$ such that
$$F\x G\sub U^n_{\le g_k(n)}\sub U^n_{\le g(n)}.\qedhere$$
\epf

We do not know whether the hypothesis in the following corollary is necessary.

\bcor[$\fb\le\cov(\cM)$]
Every $\fb$-scale set is:
\be
\itm Productively hereditarily Lindel\"of $\sone(\Omega,\gimel(\Ga))$;
\itm Productively hereditarily Lindel\"of $\sone(\Omega,\Omega^\mathrm{gp})$;
\itm Productively hereditarily Lindel\"of $\sone(\Ga,\gimel(\Ga))$; and
\itm Productively hereditarily Lindel\"of $\sone(\Ga,\gimel(\Ga))^\uparrow$.
\ee
\ecor
\bpf
The second assertion follows from the first, since
$\sone(\Omega,\Omega^\mathrm{gp})=\sone(\Omega,\gimel(\Ga))^\uparrow$ \cite{coc7}.
Similarly, the fourth assertion follows from the third.

(1) Let $X$ be a $\fb$-scale set, and let $Y$ be hereditarily Lindel\"of $\sone(\Omega,\gimel(\Ga))$.
As $X$ is $\fb$-concentrated, it is in particular $\cov(\cM)$-concentrated. By \cite[Theorem 3.1(2)]{AddGN},
$X\x Y$ satisfies $\sone(\Op,\Op)$. By Theorem \ref{bscalethm}, $X\x Y$ satisfies $\ufin(\Op,\Ga)$.
To conclude, recall that $\sone(\Omega,\gimel(\Ga))=\ufin(\Op,\Ga)\cap\sone(\Op,\Op)$ \cite{coc7}.

(3) Similar, using Proposition \ref{soneggg}, Theorem \ref{concgimel} and Theorem \ref{bscalethm}.
\epf

\brem
The only role of our restriction to hereditarily Lindel\"of in the results of this section is
to guarantee that the product with a scale set remains Lindel\"of.
\erem

\section{Combinatorial characterizations of $\ufin(\Op,\Ga)$, $\ufin(\Op,\Om)$, and $\sfin(\Op,\Op)$}\label{sec:chars}

We provide here characterizations of $\ufin(\Op,\Ga)$, $\ufin(\Op,\Om)$, and $\sfin(\Op,\Op)$ for arbitrary
topological spaces. These characterizations will be used in the following section. In this section only, the spaces
are not assumed to be Tychonoff, so that the characterizations may find additional future applications in more general contexts.
Replacing \emph{upper continuous} by \emph{continuous}
and restricting attention to separable, metrizable, zero-dimensional spaces,
the first two items in each of our characterizations become the celebrated characterizations of
Hurewicz--Rec\l{}aw \cite{Rec94} (cf.\ \cite{huremen1}).

\bdfn
Let $X$ be a topological space. For each $m$ and $n$, consider the basic open set
$$O^n_m=\set{f\in\NN}{f(n)\le m}=\pi_n\inv[\{1,\dots,m\}]$$
in $\NN$, where $\pi_n\colon\NN\to\N$ is the projection on the $n$-th coordinate.
A function $\Psi\colon X\to \NN$ is \emph{upper continuous} if, for each $n$,
the set
$$\Psi\inv[O^n_m]=\set{x\in X}{\Psi(x)(n)\le m}$$
is open in $X$.

A set-valued map $\Psi$ from $X$ to $\NN$ is \emph{principal} if there is a function $\psi\colon X\to \NN$
such that
$$\Psi(x) = \set{f\in\NN}{\forall n,\ f(n)\le \psi(x)(n)}.$$
\edfn

We use \emph{cusco} as abbreviation for \emph{compact-valued upper semicontinuous}.

\blem\label{dictionary}
Let $X$ be a topological space.
\mbox{}
\be
\itm Every cusco map from $X$ to $\NN$ is dominated by a principal one.
\itm A function $\psi\colon X\to \NN$  is upper continuous if and only if
the principal set-valued map $\Psi$ from $X$ to $\NN$ determined by $\psi$
is cusco. \qed
\ee
\elem

The equivalence of (1) and (4) in the following theorem is established in \cite[Theorem 8]{SF1}.
Our proof is, perhaps, more transparent.

\bthm\label{thm:Hchar}
Let $X$ be a topological space. The following assertions are equivalent:
\be
\itm $X$ satisfies $\ufin(\Op,\Ga)$;
\itm $X$ is Lindel\"of, and every upper continuous image of $X$ in $\NN$ is bounded;
\itm $X$ is Lindel\"of, and every principal cusco image of $X$ in $\NN$ is bounded;
\itm $X$ is Lindel\"of, and every cusco image of $X$ in $\NN$ is bounded.
\ee
\ethm
\bpf
The equivalence of $(2),(3),(4)$ follows from Lemma \ref{dictionary}.

$(1)\Impl(2)$: Let  $\Psi\colon X\to \NN$ be upper continuous.
For each $n$ and each $m$, let
$$U^n_m=\set{x\in X}{\Psi(x)(n)\le m}.$$
The sets $U^n_m$ increase with $m$, and $\set{U^n_m}{m\in\N}$
is an open cover of $X$. We may assume that $U^n_m\neq X$ for all $m$. (Otherwise,
treat the indices $m$ with $U^n_m=X$.)

Applying $\ufin(\Op,\Ga)$, there are $m_1,m_2,\dots$
such that $\set{U^n_{m_n}}{n\in\N}\in\Ga(X)$.
For each $x\in X$, $x\in U^n_{m_n}$,
and thus $\Psi(x)(n)\le m_n$, for all but finitely many $n$. In other words,
$\Psi[X]$ is bounded by the function $g(n)=m_n$.

$(2)\Impl (1)$: For each $n$, let $\cU_n=\set{U^n_m}{m\in\N}$ be an open cover
of $X$. We may assume that the sets $U^n_m$ increase with $m$.
For each $x\in X$, define
$$\Psi(x)(n)=\min\set{m}{x\in U^n_m}$$
for all $n$.
Then $\Psi$ is upper continuous.
Indeed, for each $n$,
$$\Psi\inv[\{1,\dots,m\}]=\Psi\inv[\{1\}]\cup\cdots\cup\Psi\inv[\{m\}]=
U^n_1\cup\cdots\cup U^n_m=U^n_m$$
is open in $X$.

Let $g\in\NN$ be a bound of $\Psi[X]$.
Then $\set{U^n_{g(n)}}{n\in\N}\in\Ga(X)$.
Indeed, for each $x\in X$, $\Psi(x)(n)\le g(n)$, and thus $x\in U^n_{g(n)}$, for all but finitely
many $n$.
\epf

Similarly, we have the following.

\bthm\label{thm:Schar}
Let $X$ be a topological space. The following assertions are equivalent:
\be
\itm $X$ satisfies $\ufin(\Op,\Om)$;
\itm $X$ is Lindel\"of, and no upper continuous image of $X$ in $\NN$ is finitely dominating;
\itm $X$ is Lindel\"of, and no principal cusco image of $X$ in $\NN$ is finitely dominating;
\itm $X$ is Lindel\"of, and no cusco image of $X$ in $\NN$ is finitely dominating. \qed
\ee
\ethm

\bthm\label{thm:Mchar}
Let $X$ be a topological space. The following assertions are equivalent:
\be
\itm $X$ satisfies $\sfin(\Op,\Op)$;
\itm $X$ is Lindel\"of, and no upper continuous image of $X$ in $\NN$ is dominating;
\itm $X$ is Lindel\"of, and $\NN$ is not a principal cusco image of $X$;
\itm $X$ is Lindel\"of, and $\NN$ is not a principal cusco image of $X$. \qed
\ee
\ethm

While cusco images preserve the properties mentioned above, upper continuous images need not.
However, upper continuous images are combinatorially easier to handle.

\section{Productively Lindel\"of spaces}

We conclude this paper with the following theorems concerning the property of being productively Lindel\"of.
A topological space has \emph{countable type} if each compact set in $X$ is contained in one of countable outer character.
We will use the following lemmata.

\blem[Alas--Aurichi--Junquiera--Tall \cite{AAJT11}]\label{aajt}
Let $X$ be a Lindel\"of space of countable type.
If there is an uncountable set $A\sub X$ such that $A\cap K$ is countable for every compact $K\sub X$,
then $X$ is not productively Lindel\"of.
\elem

Improving upon earlier results by several authors, Aurichi and Tall \cite{AurichiTall12}
proved that, if $\fd=\aleph_1$, then all productively Lindel\"of countable type spaces satisfy $\ufin(\Op,\Ga)$.
The following theorem both strengthens and generalizes this result.

\bthm[$\fd=\aleph_1$]
Every productively Lindel\"of metric (or just countable type) space is:
\be
\itm Productively $\ufin(\Op,\Ga)$;
\itm Productively $\sfin(\Om,\Om^\mathrm{gp})$;
\itm Productively $\sfin(\Op,\Op)$;
\itm Productively $\sfin(\Om,\Om)$; and
\itm productively $\ufin(\Op,\Om)$.
\ee
\ethm
\bpf
(2) follows from (1), and (4) from (3).

Let $X$ be a productively Lindel\"of space and let $\set{s_\alpha}{\alpha<\aleph_1}\sub\NN$ be a scale.

(1) Assume that $Y$ is $\ufin(\Op,\Ga)$ and $X\x Y$ is not.
By Theorem \ref{thm:Hchar}, there is an upper continuous $\Psi\colon X\x Y\to\NN$ such that $\Psi[X\x Y]$ is unbounded.
For each $\alpha<\aleph_1$, pick $(x_\alpha,y_\alpha)\in X\x Y$ such that
$$s_\alpha(n)\le\Psi(x_\alpha,y_\alpha)(n)$$
for infinitely many $n$.
Let $A=\set{x_\alpha}{\alpha<\aleph_1}$.
Let $K\sub X$ be compact. Then $K\x Y$ satisfies $\ufin(\Op,\Ga)$, and thus $\Psi[K\x Y]$ is bounded
by some $g\in\NN$. Let $\alpha<\aleph_1$ be such that $g<^* s_\alpha$. For each $\beta\ge\alpha$,
$$g(n)< s_\alpha(n)\le s_\beta(n)\le\Psi(x_\beta,y_\beta)(n),$$
for infinitely many $n$, and therefore $(x_\beta,y_\beta)\notin K\x Y$.
In particular, $A\cap K$ is countable.
By Lemma \ref{aajt}, $X$ is not productively Lindel\"of.

(3) Similar to (1):
Assume that $Y$ is $\sfin(\Op,\Op)$ and $X\x Y$ is not.
By Theorem \ref{thm:Mchar}, there is an upper continuous $\Psi\colon X\x Y\to\NN$ such that $\Psi[X\x Y]$ is dominating.
For each $\alpha<\aleph_1$, pick $(x_\alpha,y_\alpha)\in X\x Y$ with
$$s_\alpha\le^*\Psi(x_\alpha,y_\alpha).$$
Let $A=\set{x_\alpha}{\alpha<\aleph_1}$.
Let $K\sub X$ be compact. As $K\x Y$ satisfies $\sfin(\Op,\Op)$, $\Psi[K\x Y]$ is not dominating.
Let $g\in\NN$ be a witness for that. Let $\alpha<\aleph_1$ be such that $g<^* s_\alpha$. For each $\beta\ge\alpha$,
$$g<^* s_\alpha\le^* s_\beta\le^* \Psi(x_\beta,y_\beta),$$
and therefore $(x_\beta,y_\beta)\notin K\x Y$.
Thus, $A\cap K$ is countable, and by Lemma \ref{aajt}, $X$ is not productively Lindel\"of.

(5) Assume that $Y$ is $\ufin(\Op,\Om)$ and $X\x Y$ is not.
By Theorem \ref{thm:Schar}, there is an upper continuous $\Psi\colon X\x Y\to\NN$ such that $\Psi[X\x Y]$ is finitely dominating.
Then there is $k$ such that $\Psi[X\x Y]$ is $k$-dominating, that is,
for each $f\in\NN$ there are $(x_1,y_1),\dots,(x_k,y_k)\in X\x Y$ such that
$f\le^* \max\{\Psi(x_1,y_1),\dots,\Psi(x_k,y_k)\}$.
For each $\alpha<\aleph_1$, pick $(x^\alpha_1,y^\alpha_1),\dots,(x^\alpha_k,y^\alpha_k)\in X\x Y$ with
$$s_\alpha\le^* \max\{\Psi(x^\alpha_1,y^\alpha_1),\dots,\Psi(x^\alpha_k,y^\alpha_k)\}.$$
Define an upper continuous map
\begin{eqnarray*}
\Phi\colon X^k\x Y & \to & \NN\\
(x_1,\dots,x_k,y) & \mapsto & \max\{\Psi(x_1,y),\dots,\Psi(x_k,y)\}
\end{eqnarray*}

Let $A=\set{((x^\alpha_1,\dots,x^\alpha_k)}{\alpha<\aleph_1}\sub X^k$.
Let $K\sub X^k$ be compact.
As $K\x Y$ satisfies $\ufin(\Op,\Om)$, $\Phi[K\x Y]$ is not finitely dominating.
Let $g\in\NN$ be a witness for that. Let $\alpha<\aleph_1$ be such that $g<^* s_\alpha$. For each $\beta\ge\alpha$,
\begin{eqnarray*}
g & <^* & s_\alpha\le^* s_\beta\le^* \max\{\Psi(x^\beta_1,y^\beta_1),\dots,\Psi(x^\beta_k,y^\beta_k)\}\\
& \le^* & \max\{\Phi(x^\beta_1,\dots,x^\beta_k,y^\beta_1),\dots,\Phi(x^\beta_1,\dots,x^\beta_k,y^\beta_k)\},
\end{eqnarray*}
and therefore $(x^\beta_1,\dots,x^\beta_k)\notin K$.
Thus, $A\cap K$ is countable, and by Lemma \ref{aajt}, $X^k$ is not productively Lindel\"of.
It follows that $X$ is not productively Lindel\"of.
\epf

Say that a topological space $X$ is \emph{monotonically $\ufin(\Op,\Ga)$} if
for each upper continuous map $\Psi\colon X\to\NN$ and each $\fb$-scale $B=\set{b_\alpha}{\alpha<\fb}$, there is $\alpha<\fb$
such that for each $x\in X$,
$\Psi(x)(n)\le b_\alpha(n)$
for infinitely many $n$.
Every $\ufin(\Op,\Ga)$ space is monotonically $\ufin(\Op,\Ga)$ and every
monotonically $\ufin(\Op,\Ga)$ space is $\sfin(\Op,\Op)$.

If every $\fb$-scale is dominating (this holds, for example, in the Laver model)
then every $\sfin(\Op,\Op)$ space is monotonically $\ufin(\Op,\Ga)$.
One can prove, in ZFC, that if every $\fb$-scale is dominating, then $\aleph_1<\fb=\fd$.
On the other hand, if $\fb<\fd$ then every $\fb$-scale is an example
of a  $\sfin(\Op,\Op)$ space that is not monotonically $\ufin(\Op,\Ga)$.

\bthm[$\fb=\aleph_1$]\label{t_2}
Every productively Lindel\"of metric  (or just countable type) space is productively monotonically $\ufin(\Op,\Ga)$.
\ethm
\bpf
Let $Y$ be a monotonically $\ufin(\Op,\Ga)$ space.
Suppose, contrary to our claim, that $X\x Y$ is not monotonically $\ufin(\Op,\Ga)$.
Then there are an upper continuous map $\Psi\colon X\x Y\to\NN$ and a $\fb$-scale $B=\set{b_\alpha}{\alpha<\fb}$
such that, for each $\alpha<\aleph_1$, there is $(x_\alpha, y_\alpha)\in X\x Y$ such that
$b_\alpha\le^*\Psi(x_\alpha,y_\alpha)$.

Let $A=\set{x_\alpha}{\alpha<\aleph_1}$.
By Lemma \ref{aajt}
it is enough to prove that $A\cap K$ is countable for all compact $K\sub X$.
Indeed, it is easy to see that $\Psi$ witnesses that, for every $B$ such that $ |B\cap A|=\aleph_1$,
$B\x Y$ is not monotonically Hurewicz.
It suffices to observe that every compact space is productively monotonically  $\ufin(\Op,\Ga)$.
\epf

Alas, Aurichi, Junqueira, and Tall proved in \cite{AAJT11} that,
if $\fb=\aleph_1$, then every productively Lindel\"of  countable type space satisfies $\sfin(\Op,\Op)$.
We obtain a stronger result.

\bcor[$\fb=\aleph_1$]
Let $X$ be a productively Lindel\"of metric (or just countable type) space.
For each $\ufin(\Op,\Ga)$ space $Y$, the product space $X\x Y$ satisfies $\sfin(\Op,\Op)$.\qed
\ecor

It is known \cite{RZ12} that, if there is a Michael space, then all productively Lindel\"of spaces satisfy $\sfin(\Op,\Op)$.

\bprb
Assume that there is a Michael space. Is every productively Lindel\"of space productively $\sfin(\Op,\Op)$?
%Productively $\ufin(\Op,\Ga)$?
\eprb

\ed

\forget
I've found a simple (modulo citations) example of a sigma-compact space
and an upper continuous map in your sense
whose image is not Menger.

By Lemma 14 in our last jems paper there exists a large cover $\U$ of
the Cantor space $2^\w$ such that $\U(2^\w)$ is not Hurewicz.
Notice that since $\U(2^\w)$ is actually analytic, being not Hurewicz
means actually being not Menger.

(We are almost done, but let me complicate the situation a bit:
I am afraid of a pitfall of vacuous cases of some definitions
(like a space shouldn't be an element of its w-covers, etc.)
 if I stay with a compact space and finite covers).

Let $A={a_n:n\in\w}$ be a countable discrete space disjoint from
$2^\w$ and $X=A\cup 2^w$ be the topological sum. Then $X$ is $\sigma$-compact.
Let us write $\U$ in the form ${U_n:n\in\w}$ and let
$<k_n:n\in\w>$ be an increasing number sequence
such that $2^\w\subset\bigcup_{k_n\leq l<k_{n+1}}U_l$
for all $n\in\w$. Now set
$\V_n={V^n_l:l\in\w}$, where $V^n_l=U_{k_n+l}$
for $l<k_{n+1}-k_n$ and $V^n_l=\{ a_{l-(k_{n+1}-k_n)}  \}$
for $l\geq k_{n+1}-k_n$.

By the definition, each $\V_n$ is a cover of $X$.

Moreover, a direct verification shows that if $\Phi$ is the Marczewski
function from
$X$ to $\mathbb N$ generated  by the sequence of covers
$<\V_n:n\in\w>$, then the image $\Phi(X)$ is homeomorphic to the
topological sum of $\U(2^\w)$ and a countable discrete space,
 and hence is not Menger.
\forgotten